\documentclass{aptpub}

\authornames{Doukhan and Neumann} % insert the authors here for use in running head
\shorttitle{Mixing of GARCH-type processes} % insert short title here for use in running head

\RequirePackage[OT1]{fontenc}
\RequirePackage[authoryear]{natbib}
\RequirePackage[colorlinks,citecolor=blue,urlcolor=blue]{hyperref}

\usepackage[authoryear]{natbib}
\usepackage{color, MnSymbol, bm}
\usepackage{multirow}

\usepackage[authoryear]{natbib}
\usepackage{hyperref}
\usepackage[latin1]{inputenc}
\usepackage{graphicx}
\usepackage{bbm}
\usepackage{subcaption}
\usepackage{epstopdf}

\def\N{{\mathbb N}}
\def\Z{{\mathbb Z}}
\def\R{{\mathbb R}}
\def\P{{\mathbb P}}
\def\E{{\mathbb E}}

\def\1{{\mathbbm{1}}}
\def\wtl{{\widetilde{\lambda}}}
\def\wty{{\widetilde{Y}}}
\def\wtX{{\widetilde{X}}}

\def\esssup{\mbox{ess}\,\sup}

\let\phi\varphi
\let\epsilon\varepsilon

\newcommand{\be}{\begin{equation}}
\newcommand{\bd}{\begin{displaymath}}
\newcommand{\ed}{\end{displaymath}}
\newcommand{\bea}{\begin{eqnarray}}
\newcommand{\eea}{\end{eqnarray}}
\newcommand{\bean}{\begin{eqnarray*}}
\newcommand{\eean}{\end{eqnarray*}}

\numberwithin{equation}{section}  % If you number theorems, etc. within sections,
\begin{document}

\title{Absolute regularity of semi-contractive GARCH-type processes} % insert title - use \\ if it requires more than one line.

\authorone[Universit\'e Cergy-Pontoise]{Paul Doukhan} % Affiliation is just the name of your university or institution

\addressone{UMR 8088 Analyse, G\'eom\'etrie et Mod\'elisation\\
2, avenue Adolphe Chauvin\\
95302 Cergy-Pontoise Cedex\\
France} % Your postal address goes here.

\authortwo[Friedrich-Schiller-Universit\"at Jena]{Michael H.~Neumann}

\addresstwo{Friedrich-Schiller-Universit\"at Jena\\
Institut f\"ur Mathematik\\
Ernst-Abbe-Platz 2\\
07743 Jena\\
Germany}

\begin{abstract}
% text of abstract goes here!
We prove existence and uniqueness of a stationary distribution and absolute regularity
for nonlinear GARCH and INGARCH models of order~$(p,q)$.
In contrast to previous work we impose, besides a geometric drift condition,
only a semi-contractive condition which allows us to include models which would be ruled out
by a fully contractive condition. This results in a subgeometric rather than the more
usual geometric decay rate of the mixing coefficients. 
The proofs are heavily based on a coupling of two versions of the processes.
\end{abstract}

 % insert keywords separated by a semicolon

 % insert the primary Maths Subject Classification number in the first bracket
         % and the secondary ams number(s) in the second bracket
         % e.g. \ams{60E20}{49G03;49F10}

%%%%%%%%%%%%%%%%%%%%%%%%%%%%%%%%%%%%%%%%%%%%%%%%%%%%%%%%%%%%%%%%%%%%%%%%%%%%%%%
\section{Introduction}
\label{S1}
%%%%%%%%%%%%%%%%%%%%%%%%%%%%%%%%%%%%%%%%%%%%%%%%%%%%%%%%%%%%%%%%%%%%%%%%%%%%%%%

Conditionally heteroscedastic processes are frequently used to model the evolution
of stock prices, exchange rates and interest rates. Starting with the seminal papers
by \citet{Eng82} on autoregressive conditional heteroscedastic models (ARCH)
and \citet{Bol86} on generalized ARCH, numerous variants of these models have been
proposed for modeling financial time series; see for example \citet{FZ10} for
a detailed overview.
More recently, integer-valued GARCH models (INGARCH) which mirror the structure 
of GARCH models have been proposed for modeling time series of counts;
see for example \citet{Fok12}.

In this paper, we prove existence and uniqueness of a stationary distribution under 
a time-homogeneous dynamic. As our main result, we show absolute regularity of the 
observable process under the semi-contractive condition (\ref{1.4}) rather than a more common
fully contractive condition on the volatility function. In conjunction with standard
conditions (A1) and (A3), this results in an atypical decay rate for the coefficients of absolute regularity,
\begin{equation}
\label{RePr}
\beta_n \,=\, O(\rho^{\sqrt{n}}), \quad\mbox{ for some }\quad\rho<1.
\end{equation} 
Our technique allows to obtain this strong result even for non-stationary models with a non-homogeneous dynamic,
under uniform (in $t$) versions of our regularity conditions. 
This opens a wide range of applications for modeling real data sets.

The results hold for general GARCH processes obeying the model equations
\begin{subequations}
\begin{eqnarray}
\label{1.1a}
Y_t & = & \sigma_t \varepsilon_t, \\
\label{1.1b}
\sigma_t^2 & = & f(Y_{t-1},\ldots,Y_{t-p};\sigma_{t-1},\ldots,\sigma_{t-q}).
\end{eqnarray}
\end{subequations}
Here, $(\varepsilon_t)_t$ is a sequence of i.i.d.~random variables, where $\varepsilon_t$
is independent of all lagged random variables and $\E\varepsilon_t^2=1$.
A general INGARCH process is characterized by the model equations
\begin{subequations}
\begin{equation}
\label{1.2a}
Y_t \mid {\mathcal F}_{t-1} \,=\, Q(\lambda_t),
\end{equation}
where ${\mathcal F}_s=\sigma((Y_s,\lambda_s),(Y_{s-1},\lambda_{s-1}),\ldots)$ %}
and, analogously to the GARCH case,
\begin{equation}
\label{1.2b}
\lambda_t \,=\, f(Y_{t-1},\ldots,Y_{t-p};\lambda_{t-1},\ldots,\lambda_{t-q}).
\end{equation}
\end{subequations}
Here $\{Q(\lambda)\colon\;\; \lambda\geq 0\}$ is a family of distributions on the
non-negative integers. 
An important aspect is that such models allow for a feedback mechanism in the hidden process
which often makes a parsimonious parametrization possible.
Absolute regularity ($\beta$-mixing) with a geometric decay rate of the coefficients
of standard (linear) GARCH($p$,$q$) processes was shown in the PhD thesis of \citet{Bou98}.
Geometric $\beta$-mixing for nonlinear GARCH(1,1) specifications can be found
in \citet[proposition~5]{CC02} and \citet[Theorem~3]{FZ06}.
Properties of INGARCH processes have already been studied under a fully contractive
condition,
\begin{eqnarray}
\label{1.3}
\left| f(y_1,\ldots,y_p;\lambda_1,\ldots,\lambda_q) \,-\, f(y_1',\ldots,y_p';\lambda_1',\ldots,\lambda_q') \right| \nonumber \\
\leq\, \sum_{i=1}^p a_i \left| y_i \,-\, y_i'\right|
\,+\, \sum_{j=1}^q b_j \left| \lambda_j \,-\, \lambda_j' \right|,
\end{eqnarray}
where $y_1,\ldots,y_p,y_1',\ldots,y_p'\in\N_0=\{0,1,\ldots\}$,
$\lambda_1,\ldots,\lambda_q,\lambda_1',\ldots,\lambda_q'\geq 0$, and
$a_1,\ldots,a_p$ and $b_1,\ldots,b_q$ being non-negative constants such that $\sum_{i=1}^p a_i + \sum_{j=1}^q b_j<1$.
\citet{Neu11} showed, in the case of $p=q=1$, that condition (\ref{1.3}) implies that the bivariate process $((\lambda_t,Y_t))_t$
has a unique stationary distribution and that a stationary version of the count process $(Y_t)_t$ is absolutely regular with
mixing coefficients $\beta_n=O(\rho^n)$, for some $\rho<1$.
It was also shown that the intensity process $(\lambda_t)_t$ is not strongly mixing in general
(see Remark~3 in that paper for a simple counterexample) but ergodic.
\citet{Fra10} showed in the case of $p,q\geq 1$ that there exists
a stationary distribution. Moreover, he proved $\tau$-weak dependence as defined in \citet{WD07},
again with an exponential decay of the coefficients of weak dependence. 
Also under a fully contractive condition, \citet{FRT09} analyzed linear and nonlinear version of
INGARCH(1,1) processes. Since the verification of geometric ergodicity turned out to be unclear with 
conventional Markov chain theory, these authors proved ergodicity for a perturbed version
of the original process. As the perturbations can be chosen arbitrarily small this result could be
used to derive the asymptotic distribution of parameter estimates.

We will cover both GARCH and INGARCH models and we want to stress that we
impose a contractive condition considerably weaker than (\ref{1.3}),
\begin{equation}
\label{1.4}
\big| f(y_1,\ldots,y_p;z_1,\ldots,z_q) \,-\, f(y_1,\ldots,y_p;z_1',\ldots,z_q') \big|
\,\leq\, \sum_{i=1}^q c_i | z_i-z_i' | ,
\end{equation}  
where $c_1,\ldots,c_q$ are non-negative constants with $c_1+\cdots +c_q<1$.
This allows us to consider, for example, threshold models where the function~$f$ is specified as
\begin{equation}
\label{thresholdmodel}
f(y;\lambda) \,=\, \left\{ \begin{array}{ll} 
a \,+\, b y \,+\, c \lambda, & \quad \mbox{ if } y\in [L,U], \\
a' \,+\, b' y \,+\, c' \lambda, & \quad \mbox{ if } y\not\in [L,U].
\end{array} \right. \; 
\end{equation}
Such a specification was proposed in the framework of integer-valued time series by \cite{WMH11}.
Furthermore, our semi-contractive condition also allows us to consider functions~$f$ with
\begin{displaymath}
f(y;\lambda) \,=\, g(y) \,+\, h(\lambda)
\end{displaymath}
and with only $\mbox{Lip}(h)<1$.
Note that well-established threshold models in financial mathematics such as those
proposed for example by \citet{GJR93},
\begin{displaymath}
\sigma_t^2 \,=\, \omega \,+\, \alpha\, Y_{t-1}^2 \,+\, \beta\, Y_{t-1}^2\1_{\{Y_{t-1}<0\}} \,+\, \gamma\, \sigma_{t-1}^2 ,
\end{displaymath} 
or by \citet[page~250]{FZ10},
\begin{displaymath}
\sigma_t \,=\, \omega 
\,+\, \sum_{i=1}^p (\alpha_i^+ Y_{t-i} \1_{\{Y_{t-i}>0\}} \,-\, \alpha_i^- Y_{t-i} \1_{\{Y_{t-i}<0}\}) 
\,+\, \sum_{j=1}^q \beta_j \sigma_{t-j}
\end{displaymath}
even fulfill the fully contractive condition (\ref{1.3}).

To unify our notation, we use the expression $(\lambda_t)_t$ for the hidden process in what follows,
that is $\sigma_t^2$ will be replaced by $\lambda_t$ in case of a GARCH process.
It is worth noting at this point that, although the bivariate process $((Y_t,\lambda_t))_t$ is a Markov
chain of order $p\vee q$, the process $(Y_t)_t$ does not share this property, except for
the case $q=0$ which is not of primary interest here. 

We show as our main result that the coefficients of absolute regularity of
the observable process $(Y_t)_t$ satisfy \eqref{RePr}.
Recall that $\beta_n=\sup_{k}\beta(\mathcal{F}_{-\infty}^k,\mathcal{F}_{k+n}^\infty)$ with  
$\mathcal{F}_{k}^l=\sigma(Y_s\colon \; k\le s\le l)$ where, for any couple of $\sigma$-fields $\mathcal{A}$ and $ \mathcal{B}$:
$$
\beta(\mathcal{A},\mathcal{B})
\,=\, \sup\left\{\sum_{i=1}^\ell \sum_{j=1}^m \left|\P(A_i\cap B_j)-P(A_i)\P( B_j)\right|\right\}
$$
where the supremum is taken over partitions of $ \Omega$, $(A_i)_{1\le i\le \ell}$ and $(B_j)_{1\le i\le m}$ 
subject to  $A_i \in \mathcal{A}$ for $1\le i\le \ell$, and   $B_j \in \mathcal{B}$ for $1\le j\le m$.
This subexponential rate is quite unusual and it is a consequence of the fact that we
only impose a semi-contractive rather than a fully contractive condition.

To prove this result, we construct a coupling of two versions of the bivariate process $((Y_t,\lambda_t))_t$,
both started independently at time 0 with the stationary distribution.
These two versions, $((\wty_t,\wtl_t))_t$ and $((\wty_t',\wtl_t'))_t$, are defined on a sufficiently
rich probability space $(\widetilde{\Omega},\widetilde{\mathcal F},\widetilde{\P})$.
In the context of Markov chains, such a coupling typically leads to a coalescence of the
two versions at some random time $\tau$ and $\widetilde{\P}(\tau > n)$ then serves as an estimate of $\beta_n$.
In our case, since $(Y_t)_t$ is not a Markov chain, it can well happen that $\wty_\tau=\wty_\tau'$
at some time~$\tau$ but that afterwards these two processes diverge again. This follows from
the fact that the accompanying hidden processes $(\wtl_t)_t$ and $(\wtl_t')_t$ still can attain
different values at time~$\tau$ which means that the observable processes may diverge again with
positive probability. In view of this, we have to use $\P(\wty_m\neq \wty_m' \;\mbox{ for any } m\geq n)$
as an upper estimate for~$\beta_n$.
When the two processes reach a state with
\begin{equation}
\label{1.5}
\wty_t=\wty_t',\ldots,\wty_{t-p+1}=\wty_{t-p+1}' \quad \mbox{ and }
\quad |\wtl_t-\wtl_t'| + \cdots + |\wtl_{t-q+1}-\wtl_{t-q+1}'| \leq \rho^{\sqrt{n}},
\end{equation} 
then we have~$p$ subsequent hits and the contractive condition begins to take effect which eventually
leads to the result that both processes coalesce with a (conditional) probability exceeding $1-O(\rho^{\sqrt{n}})$.
To reach such a state with the crucial property (\ref{1.5}), the two processes need several trials,
beginning at certain stopping times $\tau_1,\tau_2,\ldots$. Because of the condition
of $|\wtl_t-\wtl_t'|+ \cdots + |\wtl_{t-q+1}-\wtl_{t-q+1}'|\leq \rho^{\sqrt{n}}$ in (\ref{1.5}),
each of these trials covers in order $\sqrt{n}$ time points.
This means, up to time~$n$ there can be in order at most $\sqrt{n}$ such trials. Such a number
of successive trials ensures that a state with (\ref{1.5}) is reached before time~$n$ with a 
probability exceeding $1-O(\rho^{\sqrt{n}})$.
This might give some insight why we obtain the unusual rate of $\rho^{\sqrt{n}}$ for the
coefficients of absolute regularity.
The desired uniqueness of the stationary law follows as a by-product of the successful coupling.
The result on absolute regularity can be extended to non-stationary GARCH-type processes;
a uniform (in~$t$) version of our semi-contractive condition will ensure this. 

The paper is organized as follows. In the next section we fix and discuss our assumptions.
Our main results are based on a coupling technique which is introduced in Subsection~\ref{SS2.1}.
To make the main ideas of our proofs easily accessible, we present the consequences of this coupling
for a simple special case in Subsection~\ref{SS2.2}.
The main results are formulated in Subsection~\ref{SS2.3}. %and an extension to non-stationary models
is briefly discussed at the end of this Subsection.
 A few applications in statistics are
mentioned in Subsection~\ref{SS2.4}.
All proofs are deferred to a final Section~\ref{S3}.

%%%%%%%%%%%%%%%%%%%%%%%%%%%%%%%%%%%%%%%%%%%%%%%%%%%%%%%%%%%%%%%%%%%%%%%%%%%%%%%
\section{Assumptions and main results}
\label{S2}
%%%%%%%%%%%%%%%%%%%%%%%%%%%%%%%%%%%%%%%%%%%%%%%%%%%%%%%%%%%%%%%%%%%%%%%%%%%%%%%

We assume that the process $(Y_t)_t$, which is defined on some probability space
$(\Omega,{\mathcal F},P)$, obeys the model equations
\begin{subequations}
\begin{eqnarray} 
\label{2.1a}
Y_t\mid {\mathcal F}_{t-1} & \sim & Q(\lambda_t), \\
\label{2.1b}
\lambda_t & = & f(Y_{t-1},\ldots,Y_{t-p};\lambda_{t-1}\ldots,\lambda_{t-q}),
\end{eqnarray}
\end{subequations}
where ${\mathcal F}_s=\sigma((Y_s,\lambda_s),(Y_{s-1},\lambda_{s-1}),\ldots)$
and $\{Q(\lambda)\colon \;\; \lambda\in [0,\infty)\}$ is some family of
univariate distributions. 
Note that assumption (\ref{2.1a}) is correctly formulated
since it follows from (\ref{2.1b}) that $\lambda_t$ is ${\mathcal F}_{t-1}$-measurable.

The canonical domain of the function~$f$ is different in the two cases
of GARCH and INGARCH models. To unify notation, we define~$f$ in both
cases on $\R^p\times [0,\infty)^q$, e.g.~by a linear interpolation in the INGARCH case.
Recall that $(\lambda_t)_t$ denotes the volatility process in the case of GARCH($p$,$q$) models ((\ref{1.1a})-(\ref{1.1b}))
and the intensity process in the INGARCH($p$,$q$) case  ((\ref{1.2a})-(\ref{1.2b})).
Here, the distribution of an observable random variable $Y_t$ conditioned on the past
is $Q(\lambda_t)$, where the parameter $\lambda_t$ itself is random,
depending on lagged variables $Y_{t-1},\ldots,Y_{t-p}$ and previous values
$\lambda_{t-1},\ldots,\lambda_{t-q}$ of the (typically hidden) accompanying process $(\lambda_t)_t$.

Possible examples we have in mind are linear or nonlinear GARCH($p$,$q$) processes,
with~$\lambda_t$ being the conditional variance of the observable variable~$Y_t$,
or integer-valued GARCH processes, where $Q(\lambda)$ is often chosen to be
a Poisson distribution with intensity parameter~$\lambda$.
Existence of a one-sided version of these processes, i.e. $t\in\N$, is guaranteed since
we can construct such processes iteratively.
We will show that there exists a stationary distribution which implies by Kolmogorov's
extension theorem (see e.g. \citet{Dur91}) that also a stationary two-sided version, i.e. $t\in\Z$, does exist.
In the proof of our main result, we also use some Markov chain techniques.
The process $(Z_t)_t$ with
$Z_t=(Y_t^2,\ldots,Y_{t-p+1}^2,\sigma_t^2,\ldots,\sigma_{t-q+1}^2)$ for a GARCH($p$,$q$) model
obeying (\ref{1.1a}) and (\ref{1.1b}) as well as
$Z_t=(Y_t,\ldots,Y_{t-p+1},\lambda_t,\ldots,\lambda_{t-q+1})$ in the INGARCH($p$,$q$) case
according to (\ref{1.2a}) and (\ref{1.2b}) has this property.
In the following it turns out to be convenient to drop the first component of the random vector~$Z_t$
and we also define $X_t=(Y_{t-1}^2,\ldots,Y_{t-p+1}^2,\sigma_t^2,\ldots,\sigma_{t-q+1}^2)$
as well as $X_t=(Y_{t-1},\ldots,Y_{t-p+1},\lambda_t,\ldots,\lambda_{t-q+1})$, respectively.

We impose the following conditions:
\begin{itemize}
\item[{\bf (A1)}] ({\em Geometric drift condition})\\
There exist positive constants $a_1,\ldots,a_{p-1}$, $b_0,\ldots,b_{q-1}$, $\kappa<1$ and $a_0<\infty$ such that,
for $V((y_1,\ldots,y_{p-1};\lambda_0,\ldots,\lambda_{q-1})) \,=\,
\sum_{i=1}^{p-1} a_i y_i + \sum_{j=0}^{q-1} b_j\lambda_j$, the condition
\begin{displaymath}
\E\left( V(X_t) \mid X_{t-1} \right)
\,\leq\, \kappa \; V(X_{t-1}) \,+\, a_0
\end{displaymath}
is fulfilled with probability~1.
\item[{\bf (A2)}] ({\em Semi-contractive condition})\\
The function $f$ is measurable and there exist non-negative constants
$c_1,\ldots,c_q$ with $c_1+\cdots +c_q<1$ such that
\begin{displaymath}
|f(y_1,\ldots,y_p; \lambda_1,\ldots,\lambda_q) \,-\, f(y_1,\ldots,y_p; \lambda_1',\ldots,\lambda_q')|
\,\leq\, \sum_{i=1}^q c_i |\lambda_i-\lambda_i'|
\end{displaymath}
for all $y_1,\ldots,y_p\in\R$, $\lambda_1,\ldots,\lambda_q,\lambda_1',\ldots,\lambda_q'\geq 0$.
\item[{\bf (A3)}] ({\em Similarity condition})\\
There exists some constant $\delta\in (0,\infty)$ such that
\begin{displaymath}
\mbox{TV}(Q(\lambda),Q(\lambda')) \,\leq\, 1 \,-\, e^{-\delta|\lambda-\lambda'|} \qquad \forall \lambda,\lambda'\geq 0,
\end{displaymath}
where $\mbox{TV}(Q_1,Q_2)=\sup_{A\in{\mathcal B}}|Q_1(A)-Q_2(A)|$ denotes the total variation distance
between probability measures~$Q_1$ and~$Q_2$.
\end{itemize}

{\rem
\label{R2.1}
In the case of $p=q=1$, $X_t$ reduces to $\lambda_t$.
Condition (A1) follows from the following drift condition which is frequently used
in the context of linear and nonlinear GARCH-type models; see e.g.
\citet{Lindner} and \citet{Fra10}.
\begin{itemize}
\item[{\bf (A1')}] There exist constants $\bar{a}_0\in [0,\infty)$, and 
$\bar{a}_1,\ldots,\bar{a}_p$, $\bar{b}_1,\ldots,\bar{b}_q\in [0,1)$, with\\
$\sum_{i=1}^p \bar{a}_i + \sum_{j=1}^q \bar{b}_j<1$ such that
\begin{itemize}
\item in the GARCH($p$,$q$) case,
\begin{displaymath}
\sigma_t^2 \,\leq\, \bar{a}_0 \,+\, \bar{a}_1 Y_{t-1}^2 \,+\, \cdots \,+\, \bar{a}_p Y_{t-p}^2 \,+\, 
\bar{b}_1 \sigma_{t-1}^2 \,+\, \cdots \,+\, \bar{b}_q \sigma_{t-q}^2,
\end{displaymath}
\item in the INGARCH($p$,$q$) case,
\begin{displaymath}
\lambda_t \,\leq\, \bar{a}_0 \,+\, \bar{a}_1 Y_{t-1} \,+\, \cdots \,+\, \bar{a}_p Y_{t-p} \,+\, 
\bar{b}_1 \lambda_{t-1} \,+\, \cdots \,+\, \bar{b}_q \lambda_{t-q}.
\end{displaymath}
\end{itemize}
\end{itemize}
}
\medskip

{\rem
\label{R2.2}
Condition (A2) is the essential difference to the fully contractive condition 
imposed e.g.~in \citet{Neu11} and \citet{Tru18}.
Here, we only assume Lipschitz continuity of~$f$ w.r.t~lagged values $\lambda_{t-1},\ldots\lambda_{t-q}$.
This includes the case of threshold models where the thresholds are set on the lagged variables
of the observable process, $Y_{t-1}^2,\ldots,Y_{t-p}^2$ or $Y_{t-1},\ldots,Y_{t-p}$, respectively.
}
\medskip

{\rem
\label{R2.3}
{} With the standard specification for GARCH models, we have that
\begin{displaymath}
Y_t \mid {\mathcal F}_{t-1} \,=\, {\mathcal N}(0,\lambda_t),
\end{displaymath}
that is, $\lambda_t$ takes the role of the conditional volatility $\sigma_t^2$.
Let $p_\lambda$ be the density of a ${\mathcal N}(0,\lambda)$ distribution.
If the volatilities satisfy $\lambda_t\geq \omega$, 
then we obtain, for $0<\omega\leq \lambda\leq \lambda'$,
\begin{displaymath}
1 \,-\, \mbox{TV}({\mathcal N}(0,\lambda),{\mathcal N}(0,\lambda'))
\,=\, \int p_\lambda \wedge p_{\lambda'} \,\geq\, \sqrt{\frac{\lambda}{\lambda'}}
\,\geq\, \frac{\lambda}{\lambda'} \,\geq\, e^{-|\lambda-\lambda'|/\lambda} \,\geq\, e^{-|\lambda-\lambda'|/\omega},
\end{displaymath}
that is, the similarity condition (A3) is fulfilled with $\delta=1/\omega$
(In order to prove the third inequality in the above display, note that $1+u\leq e^{u}$, $\forall u\geq 0$,
which implies that $\lambda'/\lambda=1+(\lambda'-\lambda)/\lambda\leq e^{|\lambda'-\lambda|/\lambda}$).

While a normal distribution seems to be the dominating choice for the distribution of the
innovations in GARCH models, there exist quite a few proposals for their integer-valued
counterparts, the INGARCH models. For the sake of an easy description, let
$({\mathcal P}_t(\lambda))_{\lambda\geq 0}$, $t\in\Z$, be a sequence of independent
standard Poisson processes.
\begin{enumerate}
\item {\em Poisson seed.} \quad
If $Q(\lambda)=\mbox{Poisson}(\lambda)$, then $Y_t$ can be expressed as
$Y_t={\mathcal P}_t(\lambda_t)$.
\item {\em Mixed Poisson seed.} \quad
Here we have the specification $Y_t={\mathcal P}_t(\lambda_t Z_t)$, where $Z_t$ is a non-negative 
random variable. The special case of a Bernoulli distributed random variable~$Z_t$,
 leads to the so-called zero-inflated Poisson model  in \citet{Lam92}; it takes into account additional unobserved data.
\item {\em Compound Poisson seed.} \quad
Let $(Z_{t,i})_{t,i\geq 0}$ be a double sequence of i.i.d.~non-negative random variables.
In this case, $Y_t$ is given by $Y_t=\sum_{i=1}^{{\mathcal P}_t(\lambda_t)} Z_{t,i}$.
This process is integer-valued if $\P(Z_{t,i}\in\N_0)=1$.
\end{enumerate}
In cases 1 and 3, the similarity assumption (A3) if fulfilled with $\delta=1$; see \citet{AJ06}.
Regarding case~2, 
let $Q_{MP}(\lambda)$ denote the mixed Poisson distribution with intensity parameter~$\lambda$.
Then,
\begin{displaymath}
\mbox{TV}(Q_{MP}(\lambda),Q_{MP}(\lambda')) \,\leq\, \E\left( 1 \,-\, e^{-Z_t|\lambda-\lambda'|} \right)
\,\leq\, 1 \,-\, e^{-\delta|\lambda-\lambda'|},
\end{displaymath}
where $\delta=\E Z_t$.
}
\medskip

{\rem
\label{R2.4}
For two probability measures~$Q_1$ and~$Q_2$ on ${\mathcal B}$, let
$d_1=dQ_1/d(Q_1+Q_2)$ and $d_2=dQ_2/d(Q_1+Q_2)$ be the respective densities
w.r.t.~the dominating measure~$Q_1+Q_2$. Then
\begin{equation}
\label{tv}
\Delta \,:=\, \int d_1\wedge d_2 \, d(Q_1+Q_2) \,=\, 1 \,-\, \mbox{TV}(Q_1, Q_2).
\end{equation}
Furthermore, using the method of maximal coupling as described for example in \citet[page~15]{dH12}
we can construct, with the aid of an additional randomization, random variables~$X_1$ and~$X_2$ such that
\begin{itemize}
\item $X_1 \sim Q_1, \qquad X_2 \sim Q_2$,
\item $P( X_1=X_2 ) \,=\, \Delta$.
\end{itemize} 
Indeed, let $U$ be a random variable with a uniform distribution on~$[0,1]$.
If $U\leq \Delta$, then we choose
\begin{displaymath}
X_1 \,=\, X_2 \,=\, F^{-1}(U),
\end{displaymath}
where $F(x)\,=\, \int_{(-\infty,x]} d_1\wedge d_2 \, d(Q_1+Q_2)$.
Here and below, $H^{-1}$ denotes the generalized inverse of a generic distribution function $H$,
that is, $H^{-1}(t)=\inf\{x\colon \; H(x)\geq t\}$ (This function is sometimes denoted by $H^{\leftarrow}$). This definition makes
sense no matter if the distribution $H$ is a continuous or discrete one. 
If $U>\Delta$, then we set
\begin{displaymath}
X_1 \,=\, G_1^{-1}(U-\Delta), \qquad X_2 \,=\, G_2^{-1}(U-\Delta),
\end{displaymath}
where $G_i(x) \,=\, \int_{(-\infty,x]} (d_i - d_1\wedge d_2) \, d(Q_1+Q_2)$,
for $i=1,2$.
}
\medskip

%%%%%%%%%%%%%%%%%%%%%%%%%%%%%%%%%%%%%%%%%%%%%%%%%%%%%%%%%%%%%%%%%%%%%%%%%%%%%%%
\subsection{Definition of the coupling}
\label{SS2.1}
%%%%%%%%%%%%%%%%%%%%%%%%%%%%%%%%%%%%%%%%%%%%%%%%%%%%%%%%%%%%%%%%%%%%%%%%%%%%%%%

We use a coupling approach to prove stationarity and absolute regularity of
the GARCH-type  process.
In the case of a stationary Markov chain $(Z_t)_{t\in\N_0}$ defined on some probability space
$(\Omega,{\mathcal F},\P)$, one usually constructs, on an appropriate probability space
$(\widetilde{\Omega},\widetilde{{\mathcal F}},\widetilde{\P})$, two versions
$(\widetilde{Z}_t)_{t\in \N_0}$ and $(\widetilde{Z}'_t)_{t\in \N_0}$  of this chain
which are started at $t=0$ independently,
both with their stationary distribution. If one succeeds to construct a coupling such that
$\widetilde{\P}(\widetilde{Z}_m \neq \widetilde{Z}'_m \mbox{ for any } m\geq n )$ tends to zero as $n\to\infty$,
then the inequality
\begin{equation}
\label{2.1.1}
\beta_n \,\leq\, \widetilde{\P}(\widetilde{Z}_m \neq \widetilde{Z}'_m \mbox{ for any } m\geq n )
\end{equation}
provides an upper bound for the mixing coefficient.
However, since a Markov process in discrete time is always strongly Markovian, it actually suffices to
derive an upper estimate for $\widetilde{\P}(\widetilde{Z}_n \neq \widetilde{Z}'_n)$ 
and we can conclude that the original process $(Z_t)_{t\in\N_0}$ on
$(\Omega,{\mathcal F},\P)$ is absolutely regular with coefficients satisfying
$\beta_n\leq \widetilde{\P}(\widetilde{Z}_n \neq \widetilde{Z}'_n)$.
In our case, the process $(Y_t)_{t}$ is not a Markov chain.
Once we have constructed a coupling of $((\widetilde{Y}_t,\widetilde{\lambda}_t))_t$ and
$((\widetilde{Y}'_t,\widetilde{\lambda}'_t))_t$, we have to stick to the estimate (\ref{2.1.1}).
(Even if $\widetilde{Y}_n=\widetilde{Y}_n'$ it could well happen that
$\widetilde{\lambda}_n\neq \widetilde{\lambda}_n'$ which means that we cannot achieve
$\widetilde{Y}_{n+1}=\widetilde{Y}_{n+1}'$ with a conditional probability of~1.)
This means that we are required to find a construction where the two versions hit at some time and
stay together afterwards (they coalesce). 

Suppose that pre-sample values $\wty_0,\ldots,\wty_{1-p}$, $\wtl_0,\ldots,\wtl_{1-q}$ 
and $\wty_0',\ldots,\wty_{1-p}'$, $\wtl_0',\ldots,\wtl_{1-q}'$ are given.
The values of $\wtl_1$ and $\wtl_1'$ arise as a result of the model equation (\ref{2.1b}),
\begin{displaymath}
\wtl_1 \,=\, f(\wty_0,\ldots,\wty_{1-p};\wtl_0,\ldots,\wtl_{1-q}), \qquad
\wtl_1' \,=\, f(\wty_0',\ldots,\wty_{1-p}';\wtl_0',\ldots,\wtl_{1-q}').
\end{displaymath}
Note that the conditional distribution of $\wty_1$ given the past has to be $Q(\wtl_1)$
and that of $\wty_1'$  $Q(\wtl_1')$.
We couple the two Markov chains in such a way that $\wty_t=\wty_t'$ with a maximum conditional probability.
According to Remark~\ref{R2.4} above, we utilize a sequence $(U_t)_{t\in\N}$ 
of i.i.d.~random variables with a uniform distribution on the interval $[0,1]$, also independent of
$(\wty_0,\wty_0',\wtl_0,\wtl_0'), (\wty_{-1},\wty_{-1}',\wtl_{-1},\wtl_{-1}'),\ldots\,$.

Let $q_1=dQ(\wtl_1)/d(Q(\wtl_1)+Q(\wtl_1'))$, $q_1'=dQ(\wtl_1')/d(Q(\wtl_1)+Q(\wtl_1'))$
and $\bar{q}_1=\int q_1\wedge q_1' \, d(Q(\wtl_1)+Q(\wtl_1'))$.

If $U_1\leq \bar{q}_1$ then we define
\begin{displaymath}
\wty_1 \,=\, \wty_1' \,=\, F_1^{-1}(U_1),
\end{displaymath}
where $F_1(x)= \int_{(-\infty, x]} q_1\wedge q_1' \, d(Q(\wtl_1)+Q(\wtl_1'))$.
If $U_1> \bar{q}_1$ then we set
\begin{displaymath}
\wty_1 \,=\, G_1^{-1}(U_1-\bar{q}_1) \qquad \mbox{ and } \qquad \wty_1' \,=\, {G_1'}^{-1}(U_1-\bar{q}_1),
\end{displaymath}
where $$G_1(x)= \int_{-\infty}^x (q_1-q_1\wedge q_1') \, d(Q(\wtl_1)+Q(\wtl_1')),\
G_1'(x)= \int_{-\infty}^x (q_1'-q_1\wedge q_1') \, d(Q(\wtl_1)+Q(\wtl_1')).$$
We iterate this process in the same way.

Let $q_t=dQ(\wtl_t)/d(Q(\wtl_t)+Q(\wtl_t'))$, $q_t'=dQ(\wtl_t')/d(Q(\wtl_t)+Q(\wtl_t'))$
and $\bar{q}_t=\int q_t\wedge q_t' \, d(Q(\wtl_t)+Q(\wtl_t'))$.
Furthermore, denote by~$F_t$, $G_t$ and $G_t'$ the distribution functions of the densities
$(q_t\wedge q_t')$, $(q_t-(q_t\wedge q_t'))$ and $(q_t'-(q_t\wedge q_t'))$,
respectively. On the basis of given values 
$\wty_{t-1},\ldots,\wty_{t-p},\wtl_{t-1},\ldots,\wtl_{t-q}$ 
and $\wty_{t-1}',\ldots,\wty_{t-p}',\wtl_{t-1}',\ldots,\wtl_{t-q}'$
we set
\begin{displaymath}
\wtl_t \,=\, f(\wty_{t-1},\ldots,\wty_{t-p};\wtl_{t-1},\ldots,\wtl_{t-q}), \qquad
\wtl_t' \,=\, f(\wty_{t-1}',\ldots,\wty_{t-p}';\wtl_{t-1}',\ldots,\wtl_{t-q}'),
\end{displaymath}
as well as
\begin{subequations}
\begin{equation}
\label{eq21.1a}
\wty_t \,=\, \wty_t' \,=\, F_t^{-1}(U_t) \qquad \qquad \mbox{ if } U_t\leq \bar{q}_t
\end{equation}
and
\begin{equation}
\label{eq21.1b}
\wty_t \,=\, G_t^{-1}(U_t-\bar{q}_t), \qquad \wty_t' \,=\, {G_t'}^{-1}(U_t-\bar{q}_t)
\qquad \qquad \mbox{ if } U_t> \bar{q}_t.
\end{equation}
\end{subequations}

%%%%%%%%%%%%%%%%%%%%%%%%%%%%%%%%%%%%%%%%%%%%%%%%%%%%%%%%%%%%%%%%%%%%%%%%%%%%%%%
\subsection{A first glimpse at the consequences of the coupling}
\label{SS2.2}
%%%%%%%%%%%%%%%%%%%%%%%%%%%%%%%%%%%%%%%%%%%%%%%%%%%%%%%%%%%%%%%%%%%%%%%%%%%%%%%

To communicate the main ideas involved in the proofs in a transparent way, we first consider
the special case of an \mbox{INGARCH(1,1)} process and present a sketch of the major steps in the
proofs of the results.
For definiteness we assume that $Y_t\mid {\mathcal F}_{t-1} \sim \mbox{Poisson}(\lambda_t)$.

Note that $\mbox{TV}(\mbox{Poisson}(\lambda),\mbox{Poisson}(\lambda'))\leq 1-e^{-|\lambda-\lambda'|}$.
To see this, assume w.l.o.g.~$\lambda\leq \lambda'$.
If $Y\sim \mbox{Poisson}(\lambda)$ and $W\sim \mbox{Poisson}(\lambda'-\lambda)$ are independent,
then $Y'=Y+W\sim \mbox{Poisson}(\lambda')$.
It follows that $P(Y\neq Y')=P(W=0)=1-e^{-|\lambda-\lambda'|}$, which implies that the similarity
condition (A3) is satisfied with $\delta=1$.

Let ${\mathcal G}_t=\sigma((\wty_t,\wty_t',\wtl_t,\wtl_t'),(\wty_{t-1},\wty_{t-1}',\wtl_{t-1},\wtl_{t-1}'),\ldots)$
denote the $\sigma$-field of the $t$-past of both versions of the processes.
Suppose that $\tau$ is some stopping time and that, for some reason, $|\wtl_{\tau+1}-\wtl_{\tau+1}'|\leq K$.
Note that $\wtl_{\tau+1}$ and $\wtl_{\tau+1}'$ are both ${\mathcal G}_\tau$-measurable, where
$${\mathcal G}_\tau=\{ G\in \bigcup_{n\in\N_0} {\mathcal G}_n\colon \quad G\cap\{\tau\leq n\}\in
{\mathcal G}_n \;\; \forall n\in\N_0\}.$$
Then, according to the maximal coupling explained above,
\begin{displaymath}
\widetilde{\P}\left( \wty_{\tau+1}=\wty_{\tau+1}' \mid {\mathcal G}_\tau \right)
\,\geq\, e^{-|\wtl_{\tau+1}-\wtl_{\tau+1}'|} \,\geq\, e^{-K}.
\end{displaymath}
If in addition $\wty_{\tau+1}=\wty_{\tau+1}'$, then the contractive condition (A2) implies that
\begin{displaymath}
| \wtl_{\tau+2} \,-\, \wtl_{\tau+2}' | \,\leq\, c_1 \; K.
\end{displaymath}
Therefore, we obtain for the next step that
\begin{displaymath}
\widetilde{\P}\left( \wty_{\tau+2}=\wty_{\tau+2}' \mid \wty_{\tau+1}=\wty_{\tau+1}', {\mathcal G}_\tau \right)
\,\geq\, e^{-c_1 K} 
\end{displaymath}
and, if additionally $\wty_{\tau+2}=\wty_{\tau+2}'$,
\begin{displaymath}
| \wtl_{\tau+3} \,-\, \wtl_{\tau+3}' | \,\leq\, c_1^2 \; K.
\end{displaymath}
Proceeding in the same way we obtain that
\begin{eqnarray}
\label{eq.1}
\lefteqn{ \widetilde{\P}\left( \wty_{\tau+1}=\wty_{\tau+1}',\ldots,\wty_{\tau+M}=\wty_{\tau+M}' \mid {\mathcal G}_\tau \right) } \nonumber \\
& = & \widetilde{\P}\left( \wty_{\tau+1}=\wty_{\tau+1}' \mid {\mathcal G}_\tau \right) 
\; \widetilde{\P}\left( \wty_{\tau+2}=\wty_{\tau+2}' \mid \wty_{\tau+1}=\wty_{\tau+1}', {\mathcal G}_\tau \right)\times \nonumber \\
& & \quad \cdots \times
\widetilde{\P}\left( \wty_{\tau+M}=\wty_{\tau+M}' \mid \wty_{\tau+1}=\wty_{\tau+1}',\ldots, \wty_{\tau+M-1}=\wty_{\tau+M-1}', {\mathcal G}_\tau \right)
\nonumber \\
& \geq & e^{-K(1+c_1+\cdots +c_1^{M-1})},
\end{eqnarray}
which leads to
\begin{eqnarray}
\label{eq.2}
\P\left( \left. \wty_{\tau+m}=\wty_{\tau+m}' ,   |\wtl_{\tau+m}-\wtl_{\tau+m}'|\leq c_1^{m-1} K,
\forall m\in\N \right| {\mathcal G}_\tau \right) \nonumber \\
\geq\, e^{-K/(1-c_1)} \,\geq\, 1 - \frac{K}{1-c_1}.
\end{eqnarray}

In what follows we sketch how (\ref{eq.2}) can be used to prove absolute regularity.
Let $\widetilde{\P}_\pi$ denote the probability where $(\wty_0,\wtl_0)$ and $(\wty_0',\wtl_0')$ are
independent and distributed with their common stationary law $\pi$.
(Its existence and uniqueness is proved in Corollary~\ref{C2.1} below.)
We define the stopping time 
\begin{displaymath}
\tau^{(n)} \,=\, \inf \{t\geq 0\colon \quad |\wtl_{t+1}-\wtl_{t+1}'| \leq C\; \rho^{n^\alpha} \},
\end{displaymath}
for some $C<\infty$ and some $\alpha>0$ whose optimal choice is explained below.
We obtain from (\ref{eq.2}) that
\begin{eqnarray}
\label{eq.3}
\beta_n & \leq & \widetilde{\P}_\pi\left( \wty_m\neq \wty_m' \quad \mbox{ for any } m\geq n \right) \nonumber \\
& \leq & \widetilde{\P}_\pi\left( \wty_m\neq \wty_m' \quad \mbox{ for any } m>\tau^{(n)} \mid {\mathcal G}_{\tau^{(n)}} \right)
\,+\, \widetilde{\P}_\pi\left( \tau^{(n)} \,\geq\, n \right) \nonumber \\
& \leq & \frac{C\; \rho^{n^\alpha}}{1-c_1} \,+\, \widetilde{\P}_\pi\left( \tau^{(n)} \,\geq\, n \right).
\end{eqnarray}
It remains to derive an upper estimate for the second term on the right-hand side of~(\ref{eq.3}).
To this end, we consider subsequent trials to achieve a state with
$|\wtl_{t}-\wtl_{t}'|\leq C_1$, for some $C_1\in (0,\infty)$,
followed by subsequent hits $\wty_{t}=\wty_{t}',\ldots,\wty_{t+d_n-1}=\wty_{t+d_n-1}'$, where $d_n=[n^\alpha]$.
We define a first stopping time as
\begin{displaymath}
\tau_1 \,=\, \inf\{t\geq 0\colon \quad \wtl_{t}+\wtl_{t}' \leq C_1\}.
\end{displaymath}
(If $\wtl_0+\wtl_0'\leq C_1$, then $\tau_1=0$. Otherwise $\tau_1$ is the first arrival
time of the process $((\wtl_t,\wtl_t'))_t$ at $A:=\{(u_1,u_2)\colon \; u_1+u_2\leq C_1\}$.)
At time $\tau_1$ we have that $|\wtl_{\tau_1}-\wtl_{\tau_1}'|\leq C_1$.
According to (\ref{eq.1}), there exists some constant $C_2>0$ such that
\begin{displaymath}
\widetilde{\P}_\pi\left( \wty_{\tau_1}=\wty_{\tau_1}',\ldots, \wty_{\tau_1+d_n-1}=\wty_{\tau_1+d_n-1}' \mid {\mathcal G}_{\tau_1-1} \right)
\,\geq\, C_2.
\end{displaymath}
After such a successful trial with~$d_n$ hits we obtain from the contractive property (A2) that
\begin{equation}
\label{eq.5}
\left| \wtl_{\tau_1+d_n} \,-\, \wtl_{\tau_1+d_n}' \right| \,\leq\, C_1 \; c_1^{d_n}.
\end{equation}
This yields that
\begin{eqnarray*}
\widetilde{\P}_\pi\left( \wty_{\tau_1+m}\neq \wty_{\tau_1+m}' \mbox{ for any } m\geq d_n \mid
\wty_{\tau_1}=\wty_{\tau_1}',\ldots,\wty_{\tau_1+d_n-1}=\wty_{\tau_1+d_n-1}', {\mathcal G}_{\tau_1-1} \right) \\
\,\leq\, \frac{C_1\; c_1^{d_n}}{1-c_1},
\end{eqnarray*}
which brings us closer to the desired result.
This means, a trial which actually leads to a favorable state with (\ref{eq.5}) covers~$d_n$ time points.
Accordingly, for $i>1$ we consider the following {\em retarded} return times as starting points for the next trials:
\begin{displaymath}
\tau_i \,=\, \inf\{ t\geq \tau_{i-1}+d_n\colon \quad \wtl_t+\wtl_t' \leq C_1\}.
\end{displaymath}

Now we are in a position to derive an upper bound for $\widetilde{\P}_\pi(\tau^{(n)}\geq n)$.
We define events
\begin{displaymath}
A_i \,=\, \left\{ \wty_{\tau_i}=\wty_{\tau_i}', \ldots, \wty_{\tau_i+d_n-1}=\wty_{\tau_i+d_n-1}' \right\}.
\end{displaymath}
Since each trial covers $d_n$ time points we cannot get more than $O(n^{1-\alpha})$ different
stopping times $\tau_i$ before time~$n$.
Let $K_n=C_3 n^{1-\alpha}$, for some $C_3>0$. It follows from Lemma~\ref{stoppingtimes} that
\begin{displaymath}
\widetilde{\P}_\pi\left( \tau_{K_n}+d_n \,\geq\, n \right)
\,\leq\, \frac1{ \eta^{n-d_n} }{ \widetilde{\E}_\pi\left( \eta^{\tau_1 + (\tau_2-\tau_1) + \cdots + (\tau_{K_n}-\tau_{K_n-1})} \right) } \\
\,=\, o\left( \rho^{n^\alpha} \right),
\end{displaymath}
for some $\eta>1$ and $\rho<1$, if $C_3$ is small enough.
Therefore, and since $\widetilde{\P}_\pi(A_1^c \cap \cdots \cap A_{K_n}^c) \leq (1-C_2)^{K_n}$,
we obtain
\begin{equation}
\label{eq.6}
\widetilde{\P}_\pi\left( \tau^{(n)} \geq n \right)
\,\leq\, \widetilde{\P}_\pi\left( \tau_{K_n}+d_n \geq n \right) \,+\, \widetilde{\P}_\pi( A_1^c \cap \cdots \cap A_{K_n}^c)
\,=\, o( \rho^{n^\alpha} ) \,+\, O(\rho^{n^{1-\alpha}}),
\end{equation}
for some $\rho<1$. The first term on the right-hand side of (\ref{eq.3}) and the second one on the right-hand side
of (\ref{eq.6}) are of the same order for the choice of $\alpha=1/2$, which gives the estimate
\begin{displaymath}
\beta_n \,=\, O( \rho^{\sqrt{n}} ).
\end{displaymath}

%%%%%%%%%%%%%%%%%%%%%%%%%%%%%%%%%%%%%%%%%%%%%%%%%%%%%%%%%%%%%%%%%%%%%%%%%%%%%%%
\subsection{Main results}
\label{SS2.3}
%%%%%%%%%%%%%%%%%%%%%%%%%%%%%%%%%%%%%%%%%%%%%%%%%%%%%%%%%%%%%%%%%%%%%%%%%%%%%%%

To prove our main results we use the coupling method described in Subsection~\ref{SS2.1}. 
Recall that $((\wty_t,\wtl_t))_t$ and $((\wty_t',\wtl_t'))_t$ denote the two versions
of the process which are coupled on a suitable probability space
$(\widetilde{\Omega},\widetilde{\mathcal F},\widetilde{\P})$ according to (\ref{eq21.1a})
and (\ref{eq21.1b}). Moreover, we remind the reader that 
${\mathcal G}_t=\sigma((\wty_t,\wty_t',\wtl_t,\wtl_t'),(\wty_{t-1},\wty_{t-1}',\wtl_{t-1},\wtl_{t-1}'),\ldots)$.
The following lemma describes the core of our coupling method.

{\lem
\label{L2.1}
Suppose that (A1) to (A3) are fulfilled
and let~$\tau$ be any stopping time such that $\wty_\tau=\wty_\tau',\ldots,\wty_{\tau-p+2}=\wty_{\tau-p+2}'$.
Then
\begin{eqnarray*}
\widetilde{\P}\left(  \wty_{\tau+m} = \wty_{\tau+m}' \forall m\in\N 
\;\; \mbox{ and } \;\; \sum_{m=1}^\infty | \wtl_{\tau+m} - \wtl_{\tau+m}' | \,\leq\, \frac{1}{1-c}
\sum_{i=1}^q |\wtl_{\tau-i+2}-\wtl_{\tau-i+2}'|\;
\Big|\; {\mathcal G}_\tau \right) \\
\,\geq\, \exp\left\{ -\frac{\delta}{1-c} \;  \sum_{i=1}^q |\wtl_{\tau-i+2}-\wtl_{\tau-i+2}'\; | \right\},
\end{eqnarray*}
where $c=c_1+\cdots +c_q$.
}
\medskip

This lemma tells us that the two processes $(\wty_t)_t$ and $(\wty_t')_t$ coalesce with a 
conditional probability greater than or equal to $\exp\{-\delta K/(1-c)\}$,
where $K=\sum_{i=1}^q |\wtl_{\tau-i+2}-\wtl_{\tau-i+2}'|$. Therefore,
in order to prove the desired decay rate for the coefficients of absolute regularity,
we show that there exists a stopping time $\tau^{(n)}$ such that
$\wty_{\tau^{(n)}}=\wty_{\tau^{(n)}}',\ldots,\wty_{\tau^{(n)}-p+2}=\wty_{\tau^{(n)}-p+2}'$,
$|\wtl_{\tau^{(n)}+1}-\wtl_{\tau^{(n)}+1}'|+\cdots +|\wtl_{\tau^{(n)}-q+2}-\wtl_{\tau^{(n)}-q+2}'|\leq \rho^{\sqrt{n}}$
and that $\widetilde{\P}(\tau^{(n)}< n) = 1-O(\rho^{\sqrt{n}})$, for some $\rho<1$.
The following main result summarizes the result of our coupling method.

{\prop
\label{P2.1}
Suppose that (A1) to (A3) are fulfilled. If
\begin{displaymath}
\widetilde{\E}\left[ V(\wtX_0) + V(\wtX_0') \right] \,<\, \infty,
\end{displaymath}
then
\begin{displaymath}
\widetilde{\P}\left( \wty_{m}=\wty_{m}' \quad \forall m\geq n \quad \mbox{ and } \quad
\sum_{m=n}^\infty | \wtl_{m} - \wtl_{m}' | \,\leq\, \frac{\rho^{\sqrt{n}}}{1-c} \right)
\,=\, 1 \,-\, O( \rho^{\sqrt{n}} ).
\end{displaymath}
}
\medskip

The following two results are immediate consequences of the main Proposition~\ref{P2.1}.

{\cor
\label{C2.1}
Suppose that (A1) to (A3) are fulfilled. Then the Markov process $(Z_t)_t$ has
a unique stationary distribution~$\pi$.
}
\medskip

{\rem
\label{R2.5}
\citet{WMH11} and \citet{DDM13} also derived properties of nonlinear INGARCH(1,1) processes which are,
as in our case here, Markov chains that are not necessarily irreducible.
\citet{WMH11} used the fact that a drift condition in conjunction with the weak Feller property of the
Markov kernel ensures the existence of a stationary distribution while its uniqueness follows
from a so-called asymptotic strong Feller property. These properties were e.g.~verified
for a Poisson threshold model with an intensity function as in (\ref{thresholdmodel}).
\citet{DDM13} extended these results to more general intensity functions, including among other examples
the log-linear Poisson autoregression model introduced by \citet{FT11}.
They also focus on the intensity process and impose the weak Feller condition directly on it.
Under an additional high-level condition on two appropriately coupled versions of the Markov chain (see
their condition (A3)) they showed that the intensity process $(\lambda_t)_t$, and as a consequence
the bivariate process $((Y_t,\lambda_t))_t$ as well, possess unique stationary distributions
and that stationary versions of the processes are ergodic.
In the case of a Poisson threshold model (\ref{thresholdmodel}) they also imposed the condition
$\max\{c,c'\}<1$ in order to ensure semi-contractivity.

Under the semicontractivity condition imposed here, we cannot derive the above mentioned Feller properties
in general. On the other hand, the coupling result stated in Proposition~\ref{P2.1} compensates 
for this failure. A metric~$d$ which resembles the coupling result is given by
\begin{multline*}
d\left( (y_1,\ldots,y_p;\lambda_1,\ldots,\lambda_q), (y_1',\ldots,y_p';\lambda_1',\ldots,\lambda_q') \right)
\\
\,=\, \1\left( (y_1,\ldots,y_p)\neq (y_1',\ldots,y_p') \right)
\,+\, \sum_{i=1}^q |\lambda_i-\lambda_i'|.
\end{multline*}
It follows for arbitrary $z\in[0,\infty)^{p+q}$ that
$P^{Z_1\mid Z_0=z'}\Rightarrow P^{Z_1\mid Z_0=z}$ as $d(z',z)\rightarrow 0$, where $\Rightarrow$ indicates
weak dependence. In other words, the weak Feller property holds w.r.t.~the metric~$d$ rather than the
more usual Euclidean norm.
As can be seen in the proof of Corollary~\ref{C2.1}, we also obtain that
\begin{displaymath}
\inf\left\{ d(\zeta_n, \zeta_n')\colon \quad \zeta_n\sim P^{Z_n\mid Z_0=z}, \zeta_n'\sim P^{Z_n\mid Z_0=z'} \right\}
\,\mathop{\longrightarrow}\limits_{n\to\infty}\, 0,
\end{displaymath}
which means that the asymptotic Feller property is also fulfilled.
\medskip

The following theorem is our main result.

{\thm
\label{T2.1}
Suppose that (A1) to (A3) are fulfilled. A stationary version of the process $(Y_t)_t$
is absolutely regular ($\beta$-mixing) with coefficients satisfying
\begin{displaymath}
\beta_n \,\leq\, C\; \rho^{\sqrt{n}} \qquad \forall n\in\N,
\end{displaymath}
for some $C<\infty$ and $\rho<1$.
}
\medskip

At this point we would like to recall that the accompanying process $(\lambda_t)_t$ is not mixing in general.
The following counter-example was already given in \citet[Remark~3]{Neu11}.
In the case of an INGARCH(1,1) process, consider the specification $f(y;\lambda)=y/2+g(\lambda)$,
where~$g$ is strictly monotone and satisfies $0<\kappa_1\leq g(\lambda)<0.5$ as well as
$|g(\lambda)-g(\lambda')|\leq \kappa_2|\lambda-\lambda'|$ for all $\lambda,\lambda'\geq 0$ and some $\kappa_2<0.5$.
Then our regularity conditions (A1) to (A3) are fulfilled.
Using the fact that $g(\lambda)\in [\kappa_1,0.5)$ we obtain from $2\lambda_t=Y_{t-1}+2g(\lambda_{t-1})$
that $Y_{t-1}=[2\lambda_t]$ and, therefore, $2g(\lambda_{t-1})=2\lambda_t-[2\lambda_t]$.
This means that we can perfectly recover~$\lambda_{t-1}$
once we know the value of~$\lambda_t$. Iterating this argument we see that we can recover from~$\lambda_t$
the complete past of the hidden process $(\lambda_t)_t$. Taking into account that the above
choice of~$f$ excludes the case that this process is purely non-random we conclude that a stationary 
version of $(\lambda_t)_t$ cannot be strongly mixing, and therefore also not be absolutely regular.

However, exploiting once more our coupling idea we can show that $\lambda_t$ can be expressed as
\begin{displaymath}
\lambda_t \,=\, g(Y_{t-1},Y_{t-2},\ldots),
\end{displaymath}
for some measurable function~$g$.
This yields ergodicity of the process $(\lambda_t)_{t\in\Z}$ and also of the bivariate
process $((Y_t,\lambda_t))_{t\in\Z}$ as stated in the following lemma.

{\thm
\label{T2.2}
Suppose that (A1) to (A3) are fulfilled. Then a stationary version of the process $((Y_t,\lambda_t))_{t\in\Z}$ 
is ergodic.
}
\medskip

Compared to absolute regularity of the process $(Y_t)_t$, the ergodicity result for the accompanying process $(\lambda_t)_t$
seems to be a bit poor. However, combined with additional structural assumptions even the property of ergodicity
might prove to be sufficient for deriving asymptotic properties of statistical procedures;
see for example \citet[Section~4]{Neu11}, \citet{LN13} and \citet{LKN15}.

{\rem
It is possible to extend our result on absolute regularity to the case of a time-varying
transition mechanism, where the function~$f$ additionally depends on time. In this case,
equation (\ref{2.1b}) has to be replaced by
\begin{equation}
\label{23.1b}
\lambda_t \,=\, f_t(Y_{t-1},\ldots,Y_{t-p};\lambda_{t-1},\ldots,\lambda_{t-q})
\end{equation}
and assumption (A2) by
\begin{itemize}
\item[{\bf (A2')}] ({\em Uniform semi-contractive condition})\\
There exist non-negative constants $c_1,\ldots,c_q$ with $c_1+\cdots +c_q<1$
such that
\begin{displaymath}
|f_t(y_1,\ldots,y_p; \lambda_1,\ldots,\lambda_q) \,-\, f_t(y_1,\ldots,y_p; \lambda_1',\ldots,\lambda_q')|
\,\leq\, \sum_{i=1}^q c_i |\lambda_i-\lambda_i'|
\end{displaymath}
for all $t\geq 0$, $y_1,\ldots,y_p\in\R$, $\lambda_1,\ldots,\lambda_q,\lambda_1',\ldots,\lambda_q'\geq 0$.
\end{itemize}
We are convinced that similar results as in our paper can be proved under these conditions
and we hope that we can report on this elsewhere.
}

\medskip

%%%%%%%%%%%%%%%%%%%%%%%%%%%%%%%%%%%%%%%%%%%%%%%%%%%%%%%%%%%%%%%%%%%%%%%%%%%%%%%
\subsection{Some applications in statistics}
\label{SS2.4}
%%%%%%%%%%%%%%%%%%%%%%%%%%%%%%%%%%%%%%%%%%%%%%%%%%%%%%%%%%%%%%%%%%%%%%%%%%%%%%%

In what follows we discuss a couple of instances where absolute regularity yields powerful
uniform limit theorems, which also indicates the relevance of the present results.
Assume that a real valued  process $(Y_t)_{t\in\Z}$ is strictly stationary and strongly mixing
with coefficients satisfying $\alpha_n\leq C\rho^{\sqrt n}$, for some $C<\infty$.
If in addition $\E g(Y_0)=0$ and $\E g^2(Y_0)\ln^2 (|g (Y_0)|\vee 1)<\infty$,
then \citet{DMR94} prove the following central limit theorem in the Skorohod space~$D[0,1]$: 
\begin{displaymath}
\frac{1}{\sqrt{n}} \sum_{j=1}^{[nu]} g(Y_i)
\,\stackrel{D[0,1]}{\longrightarrow}\,  \sigma(g)W(u), 
\end{displaymath}
where $W$ is a Brownian motion and where the series $\sigma^2(g)=\sum_{j=-\infty}^\infty \E g(Y_0)g(Y_i)$ is assumed to converge.
For the detection of changes in the mean we refer to Theorems~4.1.2 and~4.1.5 of \citet{CH}.
The same volume deals in \S \,4.4 with the detection of change points for other parameters
 involving functional central limit theorems; 
\citet{DMR95} prove a corresponding result under $\beta-$mixing. 

In the non-parametric estimation frame, the specific structure of $\beta$-mixing is also fruitful. 
\citet{vien}'s covariance inequality gives relevant bounds for the centred moments of kernel type estimators
(and more general non-parametric estimators) without imposing the existence of uniformly bounded joint densities
as this is usually done under weaker strong mixing assumptions. This inequality writes
$$
n\int_{\R^d}\var \widehat f_n(x) w(x)\,dx\le \Big(1+4\sum_{i=1}^{n-1}\beta_i\Big)\sup_{x\in \R^d} \Big\{w^2(x)\sum_{j=1}^me_i^2(x) \Big\},
$$
for projection type estimators on the vector space spanned by $\{e_1,\ldots,e_n\}$ which is
an orthonormal system of $\mathbb{L}^2(\R^d,w(x)\,dx)$. 
The standard bound of such quadratic loss has order $m/n$ under weak $\beta$-mixing assumptions.
This fact was also decisive to use model selection procedures under dependence.
\citet{bar} proposed adaptive estimation and a selection procedure 
for regression models (including autoregression)  under this $\beta$-mixing condition. 
Beyond the above mentioned covariance inequality from \citet{vien}, they used the Berbee coupling
for $\beta$-mixing sequences.

%%%%%%%%%%%%%%%%%%%%%%%%%%%%%%%%%%%%%%%%%%%%%%%%%%%%%%%%%%%%%%%%%%%%%%%%%%%%%%%%%%%%%%%%%%%%
\section{Proofs}
\label{S3}
%%%%%%%%%%%%%%%%%%%%%%%%%%%%%%%%%%%%%%%%%%%%%%%%%%%%%%%%%%%%%%%%%%%%%%%%%%%%%%%

\begin{proof}[Proof of Remark~\ref{R2.1}]
Let, for non-negative $y_1,\ldots,y_{p-1}$, $\lambda_0,\ldots,\lambda_{q-1}$ and positive
$a_1,\ldots,a_{p-1}$, $b_0,\ldots,b_{q-1}$,
$$
V\left( (y_1,\ldots,y_{p-1},\lambda_0,\ldots,\lambda_{q-1}) \right)
\,=\, \sum_{i=1}^{p-1} a_i y_i \,+\, \sum_{j=0}^{q-1} b_j \lambda_j.
$$
We consider, without loss of generality, only the case of an INGARCH($p$,$q$) process since 
the proof in the GARCH($p$,$q$) case is analogous.
Recall that $X_t=(Y_{t-1},\ldots,Y_{t-p+1},\lambda_t,\ldots,\lambda_{t-q+1})$. Then
\begin{eqnarray}
\label{pr1.1}
\lefteqn{ \E\left( V(X_t) \mid X_{t-1} \right) } \nonumber \\
& = & \E\left( a_1 Y_{t-1} \,+\, \sum_{i=2}^{p-1} a_i Y_{t-i} \,+\, b_0 \lambda_t
\,+\, \sum_{j=1}^{q-1} b_j \lambda_{t-1}
\Bigg| Y_{t-2},\ldots,Y_{t-p},\lambda_{t-1},\ldots,\lambda_{t-q} \right) \nonumber \\
& \leq & a_1 \lambda_{t-1} \nonumber %\\ & &
 {} \,+\, \sum_{i=2}^{p-1} a_i Y_{t-i} \nonumber \\
& & {} \,+\, b_0 \left( \bar{a}_0 \,+\, \bar{a}_1 \lambda_{t-1} \,+\, \sum_{i=2}^p \bar{a}_i Y_{t-i}
\,+\, \sum_{j=1}^q \bar{b}_j \lambda_{t-j} \right) 
\,+\, \sum_{j=1}^{q-1} b_j \lambda_{t-j}.
\end{eqnarray}
We are going to find positive constants $a_1,\ldots,a_{p-1},b_0,\ldots,b_{q-1}$, $\kappa<1$, and $a_0<\infty$
such that the right-hand side of (\ref{pr1.1}) is smaller than or equal to
\begin{displaymath}
a_0 \,+\, \kappa\; V(X_{t-1})
\,=\, a_0 \,+\, \kappa \left( \sum_{i=2}^p a_{i-1} Y_{t-i} \,+\, \sum_{j=1}^q b_{j-1} \lambda_{t-j} \right).
\end{displaymath}
We set, w.l.o.g., $b_0=1$ and, accordingly, $a_0=\bar{a}_0$.
Condition (A1) will be fulfilled for all possible values of the involved random variables if
\begin{subequations}
\begin{eqnarray}
\label{pr1.2a} 
a_1 \,+\, b_1 \,+\, \bar{a}_1 \,+\, \bar{b}_1 & < & 1 \\
\label{pr1.2b}
\bar{b}_j \,+\, b_j & < & b_{j-1}, \quad \mbox{ for } j=2,\ldots,q-1, \\
\label{pr1.2c}
\bar{b}_q & < & b_{q-1} \\
\label{pr1.2d}
\bar{a}_i \,+\, a_i & < & a_{i-1}, \quad \mbox{ for } i=2,\ldots,p-1, \\
\label{pr1.2e}
\bar{a}_p & < & a_{p-1},
\end{eqnarray}
\end{subequations}
where the possible choice of $\kappa$ becomes apparent at the end of the proof.

Let $\bar{a}=\sum_{i=1}^p \bar{a}_i$ and $\bar{b}=\sum_{j=1}^q \bar{b}_j$.
We choose $\varepsilon>0$ such that $\bar{a}+\bar{b}+2\varepsilon<1$
and we define
\begin{eqnarray*}
a_1 & = & \bar{a} \,-\, \bar{a}_1 \,+\, \varepsilon, \\
b_1 & = & \bar{b} \,-\, \bar{b}_1 \,+\, \varepsilon.
\end{eqnarray*}
Then (\ref{pr1.2a}) is fufilled.
Furthermore, we define recursively, for any $\delta\in (0,\varepsilon/(q-2))$,
\begin{displaymath}
b_j \,=\, b_{j-1} \,-\, \bar{b}_j \,-\, \delta, \quad \mbox{ for } j=2,\ldots,q-1,
\end{displaymath}
which implies that (\ref{pr1.2b}) holds true.
Then
\begin{displaymath}
b_{q-1} \,=\, \bar{b} \,-\, \bar{b}_1 \,-\, \cdots \,-\, \bar{b}_{q-1} \,+\, \varepsilon \,-\, (q-2)\delta
\,>\, \bar{b}_q,
\end{displaymath}
which means that (\ref{pr1.2c}) is satisfied.
Moreover, we set, for $\gamma\in(0,\varepsilon/(p-2))$,
\begin{displaymath}
a_i \,=\, a_{i-1} \,-\, \bar{a}_i \,-\, \gamma, \quad \mbox{ for } i=2,\ldots,p-1.
\end{displaymath}
Then (\ref{pr1.2d}) is fulfilled.
Finally,
\begin{displaymath}
a_{p-1} \,=\, \bar{a} \,-\, \bar{a}_1 \,-\, \cdots \,-\, \bar{a}_{p-1} \,+\, \varepsilon \,-\, (p-2) \delta
\,>\, \bar{a}_p,
\end{displaymath}
which shows that (\ref{pr1.2e}) is also satisfied.

Since all inequalities (\ref{pr1.2a}) to (\ref{pr1.2e}) are fulfilled in the strict sense
we can include a factor~$\kappa<1$
which is sufficiently close to 1 on the right-hand sides, which leaves the strict inequalities intact.
This completes the proof.
\end{proof}

\begin{proof}[Proof of Lemma~\ref{L2.1}]
Recall that $\wtl_{\tau+1}$ and $\wtl_{\tau+1}'$ are ${\mathcal G}_\tau$-measurable.
Therefore, it follows from the similarity condition (A3) and the maximal coupling scheme that
\begin{displaymath}
\widetilde{\P}\left( \wty_{\tau+1} = \wty_{\tau+1}' \mid {\mathcal G}_\tau \right)
\,\geq\, e^{-\delta |\wtl_{\tau+1}-\wtl_{\tau+1}'|}.
\end{displaymath}
If now in addition $\wty_{\tau+1}=\wty_{\tau+1}'$ then we got $p$ consecutive hits
\big($\wty_\tau=\wty_\tau',\ldots,\wty_{\tau-p+2}=\wty_{\tau-p+2}'$ was assumed\big) and the
contractive property begins to take effect, which implies that
\begin{displaymath}
|\wtl_{\tau+2}-\wtl_{\tau+2}'| \,\leq\, c_1 |\wtl_{\tau+1}-\wtl_{\tau+1}'| \,+\, \cdots \,+\,
c_q |\wtl_{\tau-q+2}-\wtl_{\tau-q+2}'|.
\end{displaymath}
Again by (A3),
\begin{displaymath}
\widetilde{\P}\left( \wty_{\tau+2} = \wty_{\tau+2}' \mid {\mathcal G}_\tau, \wty_{\tau+1} = \wty_{\tau+1}' \right)
\,\geq\, e^{-\delta |\wtl_{\tau+2}-\wtl_{\tau+2}'|}
\end{displaymath}
and, if additionally $\wty_{\tau+2} = \wty_{\tau+2}'$, then
\begin{eqnarray*}
\lefteqn{ | \wtl_{\tau+3} - \wtl_{\tau+3}' | } \\
& \leq & c_1 \; | \wtl_{\tau+2} - \wtl_{\tau+2}' | \,+\, \sum_{i=2}^q c_i \; | \wtl_{\tau+3-i} - \wtl_{\tau+3-i}' | \\
& \leq & c_1 \left( c_1 |\wtl_{\tau+1}-\wtl_{\tau+1}'| \,+\, \cdots \,+\, c_q |\wtl_{\tau-q+2}-\wtl_{\tau-q+2}'| \right)
\,+\, \sum_{i=2}^q c_i \; | \wtl_{\tau+3-i} - \wtl_{\tau+3-i}' |.
\end{eqnarray*}
Iterating these calculations we obtain for all $k\in\N$ the following general formulas.
If $\wty_{\tau-p+2}=\wty_{\tau-p+2}',\ldots,\wty_{\tau+k-1}=\wty_{\tau+k-1}'$, then
\begin{equation}
\label{pl21.1}
| \wtl_{\tau+k} - \wtl_{\tau+k}'| \,\leq\, \sum_{i=1}^q d_{k,i} \; | \wtl_{\tau-i+2} - \wtl_{\tau-i+2}' |,
\end{equation}
where $d_{1,1}=1$, $d_{1,i}=0$ if $i\geq 2$, and, for $k\geq 2$,
\begin{equation}
\label{pl21.2}
d_{k,i} \,=\, \sum_{\{l\colon \, (k+i-2)/q\leq l\leq k+i-2\}}\;
\sum_{\{(i_1,\ldots,i_l)\colon \, i_1+\cdots+i_l=k+i-2\}} c_{i_1}\times\cdots\times c_{i_l}.
\end{equation}
Therefore,
\begin{displaymath}
\widetilde{\P}\left( \wty_{\tau+k}=\wty_{\tau+k}' \mid
{\mathcal G}_\tau, \wty_{\tau+1}=\wty_{\tau+1}', \ldots, \wty_{\tau+k-1}=\wty_{\tau+k-1}' \right)
\,\geq\, e^{-\delta \sum_{i=1}^q d_{k,i} |\wtl_{\tau-i+2}-\wtl_{\tau-i+2}'|}.
\end{displaymath}
This leads to
\begin{eqnarray}
\label{pl21.3}
\lefteqn{ \widetilde{\P}\left( \wty_{\tau+1}=\wty_{\tau+1}',\ldots,\wty_{\tau+m}=\wty_{\tau+m}' \mid {\mathcal G}_\tau \right) }
\nonumber \\
& = & \P\left( \wty_{\tau+1}=\wty_{\tau+1}' \mid {\mathcal G}_\tau \right)\times 
\P\left( \wty_{\tau+2}=\wty_{\tau+2}' \mid \wty_{\tau+1}=\wty_{\tau+1}', {\mathcal G}_\tau \right)\times  \nonumber \\
& & \qquad {} \cdots\times  \P\left( \wty_{\tau+m}=\wty_{\tau+m}' \mid
\wty_{\tau+1}=\wty_{\tau+1}', \ldots, \wty_{\tau+m-1}=\wty_{\tau+m-1}', {\mathcal G}_\tau \right) \nonumber \\
& \geq & e^{-\delta \sum_{i=1}^q d_{1,i} |\wtl_{\tau-i+2}-\wtl_{\tau-i+2}'|}\times
 \cdots \times e^{-\delta \sum_{i=1}^q d_{m,i} |\wtl_{\tau-i+2}-\wtl_{\tau-i+2}'|} \nonumber \\
& = & \exp\left\{-\delta \sum_{i=1}^q D_{m,i} |\wtl_{\tau-i+2}-\wtl_{\tau-i+2}'|\right\},
\end{eqnarray}
where
\begin{displaymath} 
D_{m,i} \,:=\, \sum_{k=1}^m d_{k,i} \,\leq\, \sum_{l=0}^{m+i-2} (c_1+\cdots +c_q)^l \,\leq\, \frac{1}{1-c}.
\end{displaymath}
Since $\wty_{\tau-p+2}=\wty_{\tau-p+2}',\ldots,\wty_{\tau+m}=\wty_{\tau+m}'$ means that the 
contractive property takes effect at all time points from $\tau+1$ to $\tau+m$ we obtain that
in this case
\begin{eqnarray*}
| \wtl_{\tau+1} - \wtl_{\tau+1}' | \,+\, \cdots \,+\, | \wtl_{\tau+m+1} - \wtl_{\tau+m+1}' |
& \leq & \sum_{k=1}^{m+1} \sum_{i=1}^q d_{k,i} \; | \wtl_{\tau-i+2} - \wtl_{\tau-i+2}' | \\
& \leq & \sum_{i=1}^q D_{m+1,i} | \wtl_{\tau-i+2} - \wtl_{\tau-i+2}' |.
\end{eqnarray*}
With $m\to\infty$ we conclude that
\begin{eqnarray*} 
\label{pl2.4}
\widetilde{\P}\left( \left. \wty_{\tau+m} = \wty_{\tau+m}' \forall m\in\N 
\ \mbox{ and } \ \sum_{m=1}^\infty | \wtl_{\tau+m} - \wtl_{\tau+m}' | \,\leq\, \frac{1}{1-c}
\sum_{i=1}^q |\wtl_{\tau-i+2}-\wtl_{\tau-i+2}'|
\right| {\mathcal G}_\tau \right) \nonumber \\
\,\geq\, \exp\left\{ -\frac{\delta}{1-c} \;  \sum_{i=1}^q |\wtl_{\tau-i+2}-\wtl_{\tau-i+2}'| \right\},
\end{eqnarray*}
which proves the assertion.
\end{proof}

\bigskip

\begin{proof}[Proof of Proposition~\ref{P2.1}]
In view of the result of Lemma~\ref{L2.1}, we define a stopping time as
\begin{displaymath}
\tau^{(n)}
\,=\, \inf\{t\geq 0\colon \;\; \wty_t=\wty_t',\ldots,\wty_{t-p+2}=\wty_{t-p+2}'
\ \mbox{ and } \  |\wtl_{t+1}-\wtl_{t+1}'|+\cdots +|\wtl_{t-q+2}-\wtl_{t-q+2}'|\leq \rho^{\sqrt{n}} \},
\end{displaymath}
for some $\rho\in (0,1)$.
Recall that
\begin{eqnarray*}
{\mathcal G}_t
&=& \sigma((\wty_t,\wty_t',\wtl_t,\wtl_t'),(\wty_{t-1},\wty_{t-1}',\wtl_{t-1},\wtl_{t-1}'),\ldots)\\
&=& \sigma(\wtl_{t+1},\wtl_{t+1}',(\wty_t,\wty_t',\wtl_t,\wtl_t'),(\wty_{t-1},\wty_{t-1}',\wtl_{t-1},\wtl_{t-1}'),\ldots).
\end{eqnarray*}
It follows from Lemma~\ref{L2.1} that
\begin{eqnarray}
\label{pt1.1}
\widetilde{\P}\left( \left. \wty_{\tau^{(n)}+m} = \wty_{\tau^{(n)}+m}' \forall m\in\N 
\quad \mbox{ and } \quad \sum_{m=1}^\infty | \wtl_{\tau^{(n)}+m} - \wtl_{\tau^{(n)}+m}' |
\,\leq\, \frac{\rho^{\sqrt{n}}}{1-c} \right| {\mathcal G}_{\tau^{(n)}} \right) \nonumber \\
\,\geq\, e^{-(\delta/(1-c)) \rho^{\sqrt{n}}}
\,\geq\, 1 \,-\, \frac{\delta}{1-c} \; \rho^{\sqrt{n}}.
\end{eqnarray}
Hence, it remains to estimate $\widetilde{\P}\left( \tau^{(n)} \geq n \right)$.
To this end, we define stopping times $\tau_1,\tau_2,\ldots$
which serve as starting points of subsequent trials to reach a state with
\begin{equation}
\label{pt1.3}
\wty_t=\wty_t',\ldots,\wty_{t-p+2}=\wty_{t-p+2}'
\quad \mbox{ and } \quad |\wtl_{t+1}-\wtl_{t+1}'|+\cdots +|\wtl_{t-q+2}-\wtl_{t-q+2}'|\leq \rho^{\sqrt{n}}.
\end{equation}
Recall that $\wtX_t=(\wty_{t-1}^2,\ldots,\wty_{t-p+1}^2,\widetilde{\sigma}_t^2,\ldots,\widetilde{\sigma}_{t-q+1}^2)$,
$\wtX_t'=(\wty_{t-1}^{'2},\ldots,\wty_{t-p+1}^{'2},\widetilde{\sigma}_t^{'2},\ldots,\widetilde{\sigma}_{t-q+1}^{'2})$
in the case of a GARCH($p$,$q$) model.
Furthermore, in the INGARCH($p$,$q$) case we define these quantities as
$\wtX_t=(\wty_{t-1},\ldots,\wty_{t-p+1},\wtl_t,\ldots,\wtl_{t-q+1})$,
$\wtX_t'=(\wty_{t-1}',\ldots,\wty_{t-p+1}',\wtl_t',\ldots,\wtl_{t-q+1}')$.

Let $W_t=(V(\wtX_t)+V(\wtX_t'))/2$ and
\begin{displaymath}
\tau_1 \,=\, \inf\{t\geq 0\colon \quad W_t \,\leq\, C_1^{(0)} \},
\end{displaymath}
where $C_1^{(0)}\in (0,\infty)$ is defined in the course of the proof of Lemma~\ref{stoppingtimes} below.
Then there exists some $C_2^{(0)}>0$ such that
\begin{equation}
\label{pt1.11}
\widetilde{\P}\left( \wty_{\tau_1}=\wty_{\tau_1}' \mid {\mathcal G}_{\tau_1-1} \right) \,\geq\, C_2^{(0)}.
\end{equation}
Furthermore, it follows from (A1) that there exists some $C_1^{(1)}<\infty$ and $C_3^{(1)}>0$ such that
\begin{equation}
\label{pt1.12}
\widetilde{\P}\left( W_{\tau_1+1} \leq C_1^{(1)} \mid {\mathcal G}_{\tau_1-1}, \wty_{\tau_1}=\wty_{\tau_1}' \right)
\,\geq\, C_3^{(1)}.
\end{equation}
This, in turn, yields that there exist constants $C_2^{(1)},C_3^{(2)}>0$ and $C_1^{(2)}<\infty$
such that
\begin{equation}
\label{pt1.13}
\widetilde{\P}\left( \wty_{\tau_1+1}=\wty_{\tau_1+1}' \mid {\mathcal G}_{\tau_1-1}, \wty_{\tau_1}=\wty_{\tau_1}',
W_{\tau_1+1} \leq C_1^{(1)} \right) \,\geq\, C_2^{(1)}
\end{equation}
and
\begin{equation}
\label{pt1.14}
\widetilde{\P}\left( W_{\tau_1+2} \leq C_1^{(2)} \mid {\mathcal G}_{\tau_1-1}, \wty_{\tau_1}=\wty_{\tau_1}',
\wty_{\tau_1+1}=\wty_{\tau_1+1}', W_{\tau_1+1} \leq C_1^{(1)} \right)
\,\geq\, C_3^{(2)}.
\end{equation}
Proceeding in the same way we obtain that
\begin{eqnarray}
\label{pt1.15}
\widetilde{\P}\left(\wty_{\tau_1+p-1} = \wty_{\tau_1+p-1}' \left| \begin{array}{l}
{\mathcal G}_{\tau_1-1}, \wty_{\tau_1}=\wty_{\tau_1}',\ldots,\wty_{\tau_1+p-2}=\wty_{\tau_1+p-2}',\\
W_{\tau_1+1} \leq C_1^{(1)},\ldots, W_{\tau_1+p-1}  \leq C_1^{(p-1)}
\end{array} \right) \right. \nonumber \\
\,\geq\,  C_2^{(p-1)}.
\end{eqnarray}
This leads to
\begin{equation}
\label{pt1.16}
\widetilde{\P}\Big( \wty_{\tau_1}=\wty_{\tau_1}',\ldots,\wty_{\tau_1+p-1}=\wty_{\tau_1+p-1}'
\mid {\mathcal G}_{\tau_1-1} \Big) 
\,\geq\, C_2^{(0)}\cdots C_2^{(p-1)}\; C_3^{(1)}\cdots C_3^{(p-1)} \,=:\, C_4,
\end{equation}
that is, with a probability not smaller than $C_4>0$ we reach after~$p$ steps a state with
$\wty_{\tau_1}=\wty_{\tau_1}',\ldots,\wty_{\tau_1+p-1}=\wty_{\tau_1+p-1}'$
and $\sum_{i=1}^q b_i |\wtl_{\tau_1+p-i}-\wtl_{\tau_1+p-i}'| \leq W_{\tau_1+p-1} \leq C_1^{(p-1)}$.

Now the contractive condition begins to take effect and it follows from Lemma~\ref{L2.1}
that after $D_n-p+1:=[C_5 \sqrt{n}]$ additional hits we arrive at a state with (\ref{pt1.3}),
if $C_5$ is large enough.
This actually happens with a probability bounded away from zero.
Hence, we obtain that
\begin{displaymath}
\P\left(\left. \begin{array}{l}
 \wty_{\tau_1+D_n-1}=\wty_{\tau_1+D_n-1}',\ldots,\wty_{\tau_1+D_n-p+1}=\wty_{\tau_1+D_n-p+1}'
 \\ \mbox{ and }\qquad
\displaystyle
\sum_{i=1}^q |\wtl_{\tau_1+D_n-i+1}-\wtl_{\tau_1+D_n-i+1}'|\leq\rho^{\sqrt{n}}
\end{array}  \right| {\mathcal G}_{\tau_1-1}
\right)
\,\geq\, C_6,
\end{displaymath}
for some $C_6>0$.
This means, a trial to reach a favorable state with (\ref{pt1.3}) covers $D_n$ time points.
Accordingly, for $i>1$, we consider the following retarded return times
\begin{displaymath}
\tau_i \,=\, \inf\{t> \tau_{i-1}+D_n\colon \quad W_t \,\leq\, C_1^{(0)} \}
\end{displaymath}
Now we are in a position to derive an upper bound for $\widetilde{\P}(\tau^{(n)}\geq n)$.

We define events
\begin{displaymath}
A_i \,=\, \left\{ \wty_{\tau_i+D_n-\ell}=\wty_{\tau_i+D_n-\ell}' \mbox{ for } 1\le \ell<p,
 \mbox{ and }  \sum_{j=1}^q |\wtl_{\tau_i+D_n-j+1}-\wtl_{\tau_i+D_n-j+1}'|\leq \rho^{\sqrt{n}} \right\}.
\end{displaymath}
Let $K_n=C_7 D_n$. It follows from Lemma~\ref{stoppingtimes} that
$\widetilde{\E} \eta^{\tau_1}\leq 1+ \widetilde{\E}(\eta^{\tau_1}\mid W_0>C_1^{(0)})\leq 1+\widetilde{\E} W_0$
and 
$\widetilde{\E}(\eta^{\tau_m-\tau_{m-1}}\mid {\mathcal G}_{\tau_{m-1}-1}) \leq \eta^{D_n}C
:= \rho^{D_n} (1+(a_0+\kappa C_1^{(0)})/(1-\kappa))$, which yields
\begin{eqnarray*}
\lefteqn{ \widetilde{\E} \eta^{\tau_1+(\tau_2-\tau_1)+\cdots+(\tau_{K_n}-\tau_{K_n-1})} } \\
& = & \widetilde{\E} \left[ \eta^{\tau_1+(\tau_2-\tau_1)+\cdots+(\tau_{K_n-1}-\tau_{K_n-2})} \;
\widetilde{\E}\left( \eta^{\tau_{K_n}-\tau_{K_n-1}} \mid {\mathcal G}_{\tau_{K_n-1}-1} \right) \right] \\ 
& \leq & \eta^{D_n}\;C\; \widetilde{\E} \eta^{\tau_1+(\tau_2-\tau_1)+\cdots+(\tau_{K_n-1}-\tau_{K_n-2})} \\
& \leq & \cdots \,\leq\, \eta^{D_n(K_n-1)}\; C^{K_n-1} \; (1+\widetilde{\E} W_0).
\end{eqnarray*}
This implies that
\begin{eqnarray*}
\widetilde{\P}_\pi\left( \tau_{K_n}+D_n-1 \,\geq\, n \right)
& \leq & \frac{ \eta^{D_n(K_n-1)}\; C^{K_n-1} \; (1+\widetilde{\E}_\pi W_0) }{ \eta^{n-D_n+1} } \\
& = & O\left( \eta^{ C_7D_n^2-n-1 \; C^{C_7D_n-1} } \right)
\,=\, o\left( \rho^{\sqrt{n}} \right)
\end{eqnarray*}
if $C_7<1$ is sufficiently small. 
Therefore, and since $\widetilde{\P}(A_1^c \cap\cdots\cap A_{K_n}^c)\leq (1-C_6)^{K_n}$, we obtain that
\begin{eqnarray}
\label{pt21.21}
\widetilde{\P}( \tau^{(n)} \geq n )
& \leq & \widetilde{\P}( \tau_{K_n}+D_n-1\geq n ) \,+\, \widetilde{\P}(A_1^c \cap\cdots\cap A_{K_n}^c) \nonumber \\
& = & o( \rho^{\sqrt{n}} ) \,+\, (1-C_6)^{K_n}.
\end{eqnarray}

\end{proof}

\begin{proof}[Proof of Corollary~\ref{C2.1}]
In order to prove existence of a stationary version of $(Z_t)_t$,
it would suffice to derive this property for $(X_t)_t$, where
$X_t=(Y_{t-1},\ldots,Y_{t-p+1},\lambda_t,\ldots,\lambda_{t-q+1})$.
It follows from the drift condition (A1) that conditions (F1) and (F3), and therefore (F2) as well,
in \citet{Twe88} are fulfilled. If the Markov chain were weak Feller, i.e. for any bounded and continuous function
$\varphi\colon\;\R^{p+q-1}\rightarrow\R$ the map $x\mapsto \int \varphi(y) P^{X_1\mid X_0=x}(dy)$
were continuous, then we could conclude from Theorem~2 in \citet{Twe88} that $(X_t)_t$ has a stationary
distribution. This fact has been used e.g. in \citet{DDM13} where the weak Feller property was explicitly imposed.
The Feller property can be easily shown in case of a continuous volatility/intensity function~$f$, 
however, this might fail with a discontinuous function as they appear with certain threshold models.
We show below that the missing Feller property will be compensated by the coupling result in 
Proposition~\ref{P2.1}.

First we convert the coupling result in a convergence result for the conditional distributions
$P^{Z_n\mid X_0=x}$, where~$x$ is an arbitrarily chosen point in the range of $X_0$.
Using maximal coupling as in the proof of Proposition~\ref{P2.1} we construct two versions
of the process, $(\widetilde{Z}_t)_{t\in\N_0}$ and $(\widetilde{Z}_t')_{t\in\N_0}$,
where $\widetilde{X}_0=x$ and $\widetilde{X}_0'\sim P^{X_1\mid X_0=x}$.
We obtain that
\begin{equation}
\label{pc21.1}
\widetilde{P}\left( (\widetilde{Y}_n,\ldots,\widetilde{Y}_{n-p+1})\neq
(\widetilde{Y}_n',\ldots,\widetilde{Y}_{n-p+1}') \;\; \mbox{ or } \;\; 
\sum_{j=1}^q |\widetilde{\lambda}_{n-j+1}-\widetilde{\lambda}_{n-j+1}'|
\,>\, \frac{\rho^{\sqrt{n-q+1}}}{1-c} \right) \,=\, O\left( \rho^{\sqrt{n}} \right).
\end{equation}
Now we can construct, on a suitable probability space
$(\widetilde{\widetilde{\Omega}}, \widetilde{\widetilde{\mathcal F}}, \widetilde{\widetilde{P}})$,
a sequence of random vectors $(\zeta_n)_{n\in\N}$ such that
$\zeta_n=(\zeta_{n,1},\ldots,\zeta_{n,p},\zeta_{n,p+1},\ldots,\zeta_{p+q})^T
=(\zeta_{n,Y}^T,\zeta_{n,\lambda}^T)^T\sim P^{Z_n\mid X_0=x}$ and
\begin{equation}
\label{pc21.2}
\widetilde{\widetilde{P}}\left( \zeta_{n,Y}\neq \zeta_{n+1,Y} \;\; \mbox{ or } \;\;
\| \zeta_{n,\lambda}-\zeta_{n+1,\lambda} \|_{l_1} 
\,>\, \frac{\rho^{\sqrt{n-q+1}}}{1-c} \right) \,=\, O\left( \rho^{\sqrt{n}} \right).
\end{equation}
(Given $\zeta_1,\ldots,\zeta_n$, the vector $\zeta_{n+1}$ has to be defined
according to the conditional distribution of $\widetilde{Z}_n'$ given $\widetilde{Z}_n$.)
Since $\sum_{m=n}^\infty \rho^{\sqrt{m}}=O(\sqrt{n} \rho^{\sqrt{n}})$ we obtain from
(\ref{pc21.2}) that
\begin{equation}
\label{pc21.3}
\widetilde{\widetilde{P}}\left( \zeta_{m,Y}=\zeta_{m+1,Y} \;\; \forall m\geq n
\;\; \mbox{ and } \;\; \sum_{m=n}^\infty \|\zeta_{m,\lambda}-\zeta_{m+1,\lambda}\|_{l_1}
\leq K\sqrt{n} \rho^{\sqrt{n}} \right)
\,=\, 1 \,-\, O(\sqrt{n} \rho^{\sqrt{n}}),
\end{equation}
for some $K<\infty$.
It follows that
\begin{displaymath}
\widetilde{\widetilde{P}}\left( \bigcup_{n=1}^\infty \{\omega\colon \;\;
\zeta_{m,Y}=\zeta_{m+1,Y} \;\; \forall m\geq n\} \right) \,=\, 1,
\end{displaymath}
which means that all $\zeta_{m,Y}$ are equal for $m\geq n(\omega)$, and therefore
they are eventually equal to some random vector $\zeta_Y$.
Furthermore, since
$\zeta_{N,\lambda}=(\zeta_{N,\lambda}-\zeta_{N-1,\lambda})+\cdots
+(\zeta_{n+1,\lambda}-\zeta_{n,\lambda})+\zeta_{n,\lambda}$
we obtain that
\begin{displaymath}
\limsup_{N\to\infty} \zeta_{N,i} \,-\, \liminf_{N\to\infty} \zeta_{N,i}
\,\leq\, \sum_{m=n}^\infty |\zeta_{m+1,i}-\zeta_{m,i}|\quad \forall i=p+1,\ldots,p+q.
\end{displaymath}
Hence, it follows from (\ref{pc21.3}) that
\begin{displaymath}
\widetilde{\widetilde{P}}\left( \limsup_{N\to\infty} \zeta_{N,i} = \limsup_{N\to\infty} \zeta_{N,i}
\;\; \forall i=p+1,\ldots,p+q \right) \,=\, 1,
\end{displaymath}
which implies that $\zeta_{N,\lambda}$ converges to some random vector $\zeta_\lambda$ with
probability~1.
Let $\zeta=(\zeta_Y^T,\zeta_\lambda^T)^T$ and denote by $\pi=\widetilde{\widetilde{P}}^{\zeta}$ the distribution of $\zeta$.
Let $\varphi\colon\; \R^{p+q}\rightarrow \R$ be a bounded and uniformly continuous function.
Next we show that $\pi$ is a stationary distribution of the Markov chain $(Z_t)_t$.
Since the map $y\mapsto \int \varphi(z) P^{Z_1\mid Z_0=y}(dz)$ is continuous
in the last~$q$ arguments $y_{p+1},\ldots,y_{p+q}$ we obtain that
\begin{displaymath}
\int \left[ \int \varphi(z) P^{Z_1\mid Z_0=y}(dz) \right] \widetilde{\widetilde{P}}^{\zeta_n}(dy)
\,\mathop{\longrightarrow}\limits_{n\to\infty}\,
\int \left[ \int \varphi(z) P^{Z_1\mid Z_0=y}(dz) \right] \pi(dy),
\end{displaymath}
which yields that
\begin{eqnarray*}
\lefteqn{ \left| \int \varphi(y) \, \pi(dy) \,-\, 
\int \left[ \int \varphi(z) P^{Z_1\mid Z_0=y}(dz) \right] \pi(dy) \right| } \\
& = & \lim_{n\to\infty} \left| \int \varphi(y) \, \widetilde{\widetilde{P}}^{\zeta_n}(dy) \,-\, 
\int \left[ \int \varphi(z) P^{Z_1\mid Z_0=y}(dz) \right] \widetilde{\widetilde{P}}^{\zeta_n}(dy) \right| \\
& = & \lim_{n\to\infty} \left| \int \varphi(y) \, \widetilde{\widetilde{P}}^{\zeta_n}(dy) \,-\, 
\int \varphi(y) \, \widetilde{\widetilde{P}}^{\zeta_{n+1}}(dy) \right| \,=\, 0.
\end{eqnarray*}
Hence, $\pi$ is a stationary distribution of $(Z_t)_t$.

To show uniqueness, suppose that $\pi_1$ and $\pi_2$ are two arbitrary stationary distributions.
We start the processes to be coupled such that 
$\widetilde{Z}_0\sim \pi_1$
and
$\widetilde{Z}_0'\sim \pi_2$.
(Here, it does not matter whether or not $\widetilde{Z}_0$ and $\widetilde{Z}_0'$ are independent.)
Since both $\pi_1$ and $\pi_2$ are stationary laws we have that
\begin{equation}
\label{pc21.4}
\widetilde{Z}_n \,\sim\, \pi_1 \quad \mbox{ and } \quad \widetilde{Z}_n' \,\sim\, \pi_2 \qquad \forall n\in\N.
\end{equation}
Furthermore, it follows from the geometric drift condition (A1) that
$\widetilde{\E}\big(V(\wtX_1) + V(\wtX_1')\big)<\infty$, which implies by Proposition~\ref{P2.1} that
\begin{displaymath}
\| \widetilde{Z}_n \,-\, \widetilde{Z}_n' \| \,\stackrel{\widetilde{\P}}{\longrightarrow}\, 0,
\end{displaymath}
as $n\to\infty$. This and (\ref{pc21.4}) imply that $\pi_1=\pi_2$.
\end{proof}

\begin{proof}[Proof of Theorem~\ref{T2.1}]
Let $\pi$ denote the stationary distribution of $(Z_t)_t$ and let, for $-\infty\leq s\leq t\leq \infty$,
${\mathcal F}_{s,t}^Y=\sigma(Y_s,\ldots,Y_t)$.
We start both versions of the process at time 0 independently,
with $\widetilde{Z}_0\sim\pi$ and $\widetilde{Z}_0'\sim\pi$.
We denote by $\widetilde{\P}_\pi$ and $\widetilde{\E}_\pi$ the corresponding distribution and expectation,
respectively.
Since, by (\ref{pc23.1}) below, $\lambda_t=g(Y_{t-1},Y_{t-2},\ldots)$ we have in particular that
${\mathcal F}_{-\infty,0}^Y=\sigma(Z_0,Z_{-1},\ldots)$.
We obtain that
\begin{eqnarray*}
\beta_n 
& = & \E\left( \esssup \left\{ | P(V\mid {\mathcal F}_{-\infty,0}^Y) \,-\, P(V) |\colon
\quad V\in {\mathcal F}_{n,\infty}^Y \right\} \right) \\
& = & \E\left( \esssup \left\{ | P(V\mid Y_0,Z_0,Z_{-1},\ldots ) \,-\, P(V) |\colon
\quad V\in {\mathcal F}_{n,\infty}^Y \right\} \right) \\
& \leq & \widetilde{\E}_\pi \left( \esssup \left\{ | \widetilde{\P}_\pi( (\wty_n,\wty_{n+1},\ldots)\in A\mid {\mathcal G}_0)
\,-\, \widetilde{\P}_\pi( (\wty_n',\wty_{n+1}',\ldots)\in A\mid {\mathcal G}_0) |\colon \;\; A\in {\mathcal C} \right\} \right) \\
& \leq & \widetilde{\E}_\pi \left( \widetilde{\P}_\pi( \exists m\geq n\colon \;\; \wty_m\neq \wty_m\mid {\mathcal G}_0) \right)
\,=\, \widetilde{\P}_\pi( \exists m\geq n\colon \;\; \wty_m\neq \wty_m).
\end{eqnarray*}
Here, ${\mathcal C}$ denotes the $\sigma$-field generated by the cylinder sets.
Proposition~\ref{P2.1} yields that $\beta_n=O(\rho^{\sqrt{n}})$, as required.
\end{proof}

\begin{proof}[Proof of Theorem~\ref{T2.2}]
Let $((Y_t,\lambda_t))_{t\in\Z}$ be a stationary version of the process.
We will show that there exists a measurable function~$g\colon\;\N_0^\infty\rightarrow [0,\infty)$
such that $\lambda_t=g(Y_{t-1},Y_{t-2},\ldots)$.
To this end, we consider the same ``forward iterations'' as in the proof of Lemma~\ref{L2.1}.
We use the true values $Y_0,\ldots,Y_{1-p},\lambda_0,\ldots,\lambda_{1-q}$ as well as
$Y_0,\ldots,Y_{1-p},\lambda_0',\ldots,\lambda_{1-q}'$ with $\lambda_0'=\ldots =\lambda_{1-q}'=0$
as starting values.
Then we define, according to the model equation (\ref{2.1b}),
\begin{eqnarray*}
\lambda_1 & = & f(Y_0,\ldots,Y_{1-p};\lambda_0,\ldots,\lambda_{1-q}), \\
\lambda_1' & = & f(Y_0,\ldots,Y_{1-p};\lambda_0',\ldots,\lambda_{1-q}') 
\,=:\, g^{[1]}(Y_0,\ldots,Y_{1-p}).
\end{eqnarray*}
Iterating this scheme we obtain
\begin{eqnarray*}
\lambda_k & = & f(Y_{k-1},\ldots,Y_{k-p};\lambda_{k-1},\ldots,\lambda_{k-q}), \\
\lambda_k' & = & f(Y_{k-1},\ldots,Y_{k-p};\lambda_{k-1}',\ldots,\lambda_{k-q}') 
\,=:\, g^{[k]}(Y_{k-1},\ldots,Y_{1-p}).
\end{eqnarray*}
Note that in all steps matching values of the process $(Y_t)_t$ are used for
computing $\lambda_k$ and $\lambda_k'$, which means that the contractive property takes effect
at each step. Therefore we obtain, analogously to (\ref{pl21.1}) in the proof
of Lemma~\ref{L2.1},
\begin{displaymath}
| \lambda_k \,-\, \lambda_k' |
\,\leq\, \sum_{i=1}^q d_{k+1,i} \lambda_{1-i},
\end{displaymath}
where it follows from (\ref{pl21.2}) that $d_{k+1}\rightarrow_{k\to\infty} 0$.
By stationarity we conclude, for fixed $t\in\Z$, that
\begin{displaymath}
\E\left| \lambda_t \,-\, g^{[k]}(Y_{t-1},\ldots,Y_{t-p-k+1}) \right|
\,\leq\, \sum_{i=1}^q d_{k+1,i} \E \lambda_{t-k-i+1}
\,\mathop{\longrightarrow}\limits_{k\to\infty}\, 0,
\end{displaymath}
that is, as $k\to\infty$, $g^{[k]}(Y_{t-1},\ldots,Y_{t-p-k+1})$ converges in $L_1$ to $\lambda_t$.
By taking an appropriate subsequence we also get almost sure convergence. This means that there
exists some measurable function $g\colon\;\N_0^\infty\rightarrow [0,\infty)$ such that
\begin{equation}
\label{pc23.1}
\lambda_t \,=\, g(Y_{t-1},Y_{t-2},\ldots) \qquad \mbox{almost surely}.
\end{equation}
Since absolute regularity of the process $(Y_t)_{t\in\Z}$ implies strong mixing (see e.g.~\citet[p.~20]{Dou94})
we conclude from Remark~2.6 on page~50 in combination with Proposition~2.8 on page~51 in \citet{Bra07}
that any stationary version of this process is also ergodic.

Finally, we conclude from the representation (\ref{pc23.1})
by proposition~2.10(ii) in \citet[p.~54]{Bra07} that also
the bivariate process $((Y_t,\lambda_t))_{t\in\Z}$  is ergodic.
\end{proof}

{\lem
\label{stoppingtimes}
Suppose that (A1) is fulfilled. Then
\begin{itemize}
\item[(i)] $\widetilde{\E}( \eta^{\tau_1} \mid {\mathcal G}_{-1} )
\,\leq\, ( V(\wtX_0) \,+\, V(\wtX_0') )/2, 
\quad \mbox{ if } ( V(\wtX_0) \,+\, V(\wtX_0') )/2 > C_1^{(0)}$,\\
where $\eta=2/(1+\kappa)$ and $C_1^{(0)}=(2a_0 + 2)/(1-\kappa)$.\\
\item[(ii)] $\widetilde{\E}( \eta^{\tau_{m+1}-\tau_m} \mid {\mathcal G}_{\tau_m-1} )
\,\leq\, \rho^{D_n} \left( 1 \,+\, \frac{a_0+\kappa C_1^{(0)}}{1-\kappa} \right)$.
\end{itemize}
}

\begin{proof}[Proof of Lemma~\ref{stoppingtimes}]
We already defined $\wtX_t=(\wty_{t-1}^2,\ldots,\wty_{t-p+1}^2,\widetilde{\sigma}_t^2,\ldots,\widetilde{\sigma}_{t-q+1}^2)$, and
$\wtX_t'=(\wty_{t-1}^{'2},\ldots,\wty_{t-p+1}^{'2},\widetilde{\sigma}_t^{'2},\ldots,\widetilde{\sigma}_{t-q+1}^{'2})$
in the case of a GARCH($p$,$q$) model.
Furthermore,  in the INGARCH($p$,$q$) case we set analogously
$\wtX_t=(\wty_{t-1},\ldots,\wty_{t-p+1},\wtl_t,\ldots,\wtl_{t-q+1})$,
$\wtX_t'=(\wty_{t-1}',\ldots,\wty_{t-p+1}',\wtl_t',\ldots,\wtl_{t-q+1}')$.
Let $W_t=(V(\wtX_t)+(\wtX_t'))/2$.

Since $\wty_{t-1}\mid {\mathcal G}_{t-1}=Q(\wtl_{t-1})$ we see that
$\widetilde{\E}( \wty_{t-1}\mid {\mathcal G}_{t-1} )=\widetilde{\E}( \wty_{t-1}\mid \wtX_{t-1})$
and $\widetilde{\E}( \wtl_t\mid {\mathcal G}_{t-1} )
=\widetilde{\E}( f(\wty_{t-1},\ldots,\wty_{t-p};\wtl_{t-1},\ldots,\wtl_{t-q})\mid {\mathcal G}_{t-1} )
=\widetilde{\E}( \wtl_t\mid \wtX_{t-1} )$.

Therefore we obtain
$\widetilde{\E}(V(\wtX_t)\mid {\mathcal G}_{t-1})=\widetilde{\E}(V(\wtX_t)\mid \wtX_{t-1})$
and, analogously,
$\widetilde{\E}(V(\wtX_t')\mid {\mathcal G}_{t-1})=\widetilde{\E}(V(\wtX_t)\mid \wtX_{t-1}')$.
Hence, we obtain from the geometric drift condition (A1) that
\begin{equation}
\label{pls1.0}
\widetilde{\E}\left( W_t \mid {\mathcal G}_{t-1} \right) \,\leq\, \kappa \; W_{t-1} \,+\, a_0.
\end{equation}
This implies that
\begin{equation}
\label{pls.1}
\widetilde{\E}( W_t\mid {\mathcal G}_{t-1} ) \,\leq\, \eta^{-1} W_{t-1} \,-\, 1, \quad \mbox{ if } W_{t-1}>C_1^{(0)}
\end{equation}
and 
\begin{equation}
\label{pls.2}
\widetilde{\E}( W_t\mid {\mathcal G}_{t-1} ) \,\leq\, \kappa C_1^{(0)} \,+\, a_0, \quad \mbox{ if } W_{t-1}\leq C_1^{(0)}.
\end{equation}\

In what follows we adapt the line of arguments from \citet{NT82}, who derived similar bounds for
stopping times in the context of a Markov chain.\medskip

\noindent
{\bf Proof of (i)}

Let $W_0=x>C_1^{(0)}$. We denote by $\widetilde{\P}_x$ and $\widetilde{\E}_x$
the conditional distribution and expectation, respectively, given $W_0=x$.
It follows from (\ref{pls.1}) that
\begin{displaymath}
\widetilde{\E}_x( W_1 ) \,\leq\, \eta^{-1} x \,-\, 1,
\end{displaymath}
which implies that
\begin{equation}
\label{pls.3}
x \,-\, \eta \; \E_x\left( W_1 \right) \,\geq\, \eta.
\end{equation}
Analogously we conclude from (\ref{pls.1}) that
\begin{displaymath}
\1( W_1>C_1^{(0)} ) \; \widetilde{\E}_x( W_2 \mid W_1 ) \,\leq\, \1( W_1>C_1^{(0)} ) \; (\eta^{-1} W_1 \,-\, 1),
\end{displaymath}
which yields that
\begin{displaymath}
\1( W_1>C_1^{(0)} ) \; \left( W_1 \,-\, \eta \; \widetilde{\E}_x\left( W_2\mid W_1 \right) \right)
\,\geq\, \eta \; \1( W_1>C_1^{(0)} ).
\end{displaymath}
Multiplying both sides by $\eta$ and taking the expectation over $W_1$ under the condition $W_0=x$
we obtain
\begin{equation}
\label{pls.4}
\widetilde{\E}_x\left( \1( W_1>C_1^{(0)} ) \; \left( \eta W_1 \,-\, \eta^2 W_2 \right) \right)
\,\geq\, \eta^2 \; \widetilde{\P}_x( W_1>C_1^{(0)} ).
\end{equation}
Proceeding in the same way we conclude
\begin{eqnarray}
\label{pls.5}
\widetilde{\E}_x\left( \1( W_1>C_1^{(0)},\ldots,W_k>C_1^{(0)} ) \; \left( \eta^k W_k \,-\, \eta^{k+1} W_{k+1} \right) \right)
\nonumber \\
\geq \eta^{k+1} \; \widetilde{\P}_x( W_1>C_1^{(0)},\ldots,W_k>C_1^{(0)} ).
\end{eqnarray}\
Adding both sides of (\ref{pls.3}) to (\ref{pls.5}) we obtain that
\begin{eqnarray*}
x & \geq & \sum_{k=0}^\infty \eta^{k+1} \; \widetilde{\P}_x( W_1>C_1^{(0)},\ldots,W_k>C_1^{(0)} ) \\
& = & \sum_{k=0}^\infty \eta^{k+1} \; \widetilde{\P}_x\left( \tau_1\geq k+1 \right)
\,\geq\, \widetilde{\E}_x\left( \eta^{\tau_1} \right),
\end{eqnarray*}
as required.
\medskip

\noindent
{\bf Proof of (ii)}

Here we have to take into account that $\tau_{m+1}$ is not a usual but a retarded
return time. Recall that $X_{\tau_m}$ is ${\mathcal G}_{\tau_m-1}$-measurable.
Since $X_{\tau_m}\leq C_1^{(0)}$ we obtain from (i) that
\begin{eqnarray}
\label{pls.11}
\lefteqn{ \widetilde{\E}\left( \eta^{\tau_{m+1}-\tau_m} \mid {\mathcal G}_{\tau_m-1} \right) } \nonumber \\
& = & \eta^{D_n} \widetilde{\P}\left( W_{\tau_m+D_n} \leq C_1^{(0)} \mid {\mathcal G}_{\tau_m-1} \right) \nonumber \\
& & {} \,+\, \eta^{D_n} \int_{(C_1^{(0)},\infty)}
\widetilde{\E}( \eta^{\tau_{m+1}-(\tau_m+D_n)} \mid {\mathcal G}_{\tau_m-1}, W_{\tau_m+D_n}=x )
\; \widetilde{\P}^{W_{\tau_m+D_n}\mid {\mathcal G}_{\tau_m-1}}(dx) \nonumber \\
& \leq & \eta^{D_n} \left( 1 \,+\, \widetilde{\E}( W_{\tau_m+D_n} \mid {\mathcal G}_{\tau_m-1} ) \right)
\end{eqnarray}
Furthermore, since $W_{\tau_m}\leq C_1^{(0)}$ we obtain from (\ref{pls1.0}) that 
\begin{displaymath}
\widetilde{\E}( W_{\tau_m+1} \mid {\mathcal G}_{\tau_m-1} ) 
\,\leq\, \kappa W_{\tau_m} \,+\, a_0
\,\leq\, \kappa C_1^{(0)} \,+\, a_0,
\end{displaymath}
\begin{eqnarray*}
\widetilde{\E}( W_{\tau_m+2} \mid {\mathcal G}_{\tau_m-1} )
& = & \widetilde{\E}\left( \widetilde{\E}( W_{\tau_m+2} \mid {\mathcal G}_{\tau_m-1}, W_{\tau_m+1} )
\mid {\mathcal G}_{\tau_m-1} \right) \\
& \leq & \widetilde{\E}\left( \kappa W_{t_m+1} \,+\, a_0 \mid {\mathcal G}_{\tau_m-1} \right) \\
& \leq & 2a_0 \,+\, \kappa (\kappa C_1^{(0)} \,+\, a_0 ),
\end{eqnarray*}
and, eventually,
\begin{equation}
\label{pls.12}
\widetilde{\E}( W_{\tau_m+k} \mid {\mathcal G}_{\tau_m-1} )
\,\leq\, \frac{a_0 \,+\, \kappa C_1^{(0)}}{1-\kappa} \quad \forall k\in\N.
\end{equation}
(ii) now follows from (\ref{pls.11}) and (\ref{pls.12}).
\end{proof}

\acks
This work has been developed within the MME-DII center of excellence (ANR-11-LABEX-0023-01)
and with the help of PAI-CONICYT MEC $\mbox{N}^{o}$~80170072. The first author wishes to thank the University of Jena and
Universidad of Valparaiso for their hospitality.
The research of the second author was supported by a guest professorship of IEA at the University of Cergy-Pontoise.
We thank two anonymous referees for their comments which helped us to improve the presentation of our results.

% Place the text of your acknowledgements after the \acks command.
% \acks generates the heading "Acknowledgements".
% If you wish to make only one acknowledgement, use \ack.
% \ack generates the heading "Acknowledgement".

% Reference list
%
% References should be in the following form (or the BibTeX file
% apt.bst should be used):
%
% For a journal:
% Surname, Initial (year). Title of paper. {\em Journal title}
% {\bf Vol,} page--range.
%
% For a book:
% Surname, Initial (year). {\em Book title}. Publisher, Address.
%
% Note the following example of a reference list.

\bibliographystyle{harvard}

\end{document}